\documentclass{article}

\usepackage{amstex}

\newcommand{\ssp}{\discretionary{}{}{\,}}

\def\bewname{Proof}
\def\definname{Definition}
\def\thmname{Theorem}
\def\satzname{Proposition}
\def\lemmaname{Lemma}
\def\folgname{Corollary}
\def\bemname{Remark}
\def\bemsname{Remarks}
\def\bspname{Example}
\newcommand{\mc}{\mathcal}
\newcommand{\mf}{\mathfrak}
\newcommand{\mb}{\mathbf}

\newcommand{\A}{\mathcal{A}}
\newcommand{\B}{\mathcal{B}}
\newcommand{\T}{\mathcal{T}}
\newcommand{\U}{\mathcal{U}}

\newcommand{\comp}{\mathbb{C}}
\newcommand{\OGq}{\mc{O}(G_q)}
\newcommand{\RGq}{\mc{R}(G_q)}
\newcommand{\OGqn}{\OGq^\circ}
\newcommand{\RGqn}{\RGq^\circ}
\newcommand{\Uqg}{\mc{U}_q(\mf{g})}

\newcommand{\tUqg}{{\tilde{\mc{U}}_q(\mf{g})}}

\newcommand{\Uqnp}{\mc{U}_q(\mf{n}_+)}
\newcommand{\Uqnm}{\mc{U}_q(\mf{n}_-)}
\newcommand{\X}{\mc{X}}
\newcommand{\Y}{\mc{Y}}
\newcommand{\br}{\mb{r}}

\newcommand{\liesl}{\mathrm{sl}}
\newcommand{\lieso}{\mathrm{so}}
\newcommand{\liesp}{\mathrm{sp}}
\newcommand{\smashright}{\rhd\!\!<}

\newcommand{\fuw}{\omega}
\newcommand{\chru}{T}
\newcommand{\pchru}{T_+}

\newcounter{beh}[section]

\def\thebeh{\thesection.\arabic{beh}}

\newenvironment{bew}[1]{\vspace{1ex}\textbf{\bewname{#1.}}}
{\ \hspace*{\fill}\rule{1ex}{1ex}\vspace{1ex}}

\newenvironment{bew*}[1]{\textbf{\bewname {#1.}}}
{\ \hspace*{\fill}\rule{1ex}{1ex}\vspace{1ex}}

\newenvironment{defin*}{\refstepcounter{beh}
\textbf{\definname~\thebeh.} }{\vspace{1ex}}

\newenvironment{thm}{\refstepcounter{beh}\vspace{1ex}
\textbf{\thmname~\thebeh.} \em }{\em\vspace{1ex}}

\newenvironment{thm*}{\refstepcounter{beh}
\textbf{\thmname~\thebeh.} \em }{\em\vspace{1ex}}

\newenvironment{satz}{\refstepcounter{beh}\vspace{1ex}
\textbf{\satzname~\thebeh.} \em }{\em\vspace{1ex}}

\newenvironment{satz*}{\refstepcounter{beh}
\textbf{\satzname~\thebeh.} \em }{\em\vspace{1ex}}

\newenvironment{lemma}{\refstepcounter{beh}\vspace{1ex}
\textbf{\lemmaname~\thebeh.} \em }{\em\vspace{1ex}}

\newenvironment{lemma*}{\refstepcounter{beh}
\textbf{\lemmaname~\thebeh.} \em }{\em\vspace{1ex}}

\newenvironment{folg}{\refstepcounter{beh}\vspace{1ex}
\textbf{\folgname~\thebeh.} \em }{\em\vspace{1ex}}

\newenvironment{folg*}{\refstepcounter{beh}
\textbf{\folgname~\thebeh.} \em }{\em\vspace{1ex}}

\newenvironment{bem}{\vspace{1ex}\textbf{\bemname.} }
{\ \hspace*{\fill}\rule{1ex}{1ex}\vspace{1ex}}

\newenvironment{bem*}{\textbf{\bemname.} }
{\ \hspace*{\fill}\rule{1ex}{1ex}\vspace{1ex}}

\newenvironment{bems*}{\textbf{\bemsname.} }
{\ \hspace*{\fill}\rule{1ex}{1ex}\vspace{1ex}}

\newenvironment{bsp}{\refstepcounter{beh}\vspace{1ex}
\textbf{\bspname~\thebeh.} }
{\ \hspace*{\fill}\rule{1ex}{1ex}\vspace{1ex}}

\newenvironment{bsp*}{\refstepcounter{beh}
\textbf{\bspname~\thebeh.} }
{\ \hspace*{\fill}\rule{1ex}{1ex}\vspace{1ex}}

\newcommand{\Gammalinv}{{_{\mathrm{inv}}\Gamma}}

\newcommand{\linv}[1]{{_{\mathrm{inv}}{#1}}}

\newcommand{\adR}{\mathrm{ad}_R\,}
\newcommand{\adL}{\mathrm{ad}_L\,}
\newcommand{\AdR}{\mathrm{Ad}_R\,}

\newcommand{\cont}{^{\mathrm{c}}}
\newcommand{\contr}{\mathrm{c}}
\newcommand{\dif}{\mathrm{d}}

\newcommand{\id}{\mathrm{id}}

\author{I.~Heckenberger\footnotemark[1]\,\, and
K.~Schm{\" u}dgen\footnotemark[2]}

\title{Classification of Bicovariant Differential Calculi on
the Quantum Groups $SL_q(n+1)$ and $Sp_q(2n)$}

\date{}

\begin{document}

\maketitle

\footnotetext[1]{e-mail: heckenbe@@mathematik.uni-leipzig.de}
\footnotetext[2]{e-mail: schmuedg@@mathematik.uni-leipzig.de}

\begin{abstract}
For transcendental values of $q$ all bicovariant first order
differential calculi
on the coordinate Hopf algebras of the quantum groups
$SL_q(n+1)$ and $Sp_q(2n)$ are classified. It is shown that the
irreducible bicovariant first order calculi are determined by an
irreducible corepresentation of the quantum group and a complex
number $\zeta$ such that $\zeta^{n+1}=1$ for $SL_q(n+1)$ and
$\zeta^2=1$ for $Sp_q(2n)$.
Any bicovariant calculus is inner and its quantum Lie algebra is generated
by a central element. The main technical ingredient is a result
of the Hopf algebra $\RGqn$ for arbitrary simple Lie algebras.
\end{abstract}

AMS subject classification: 17 B37, 46 L87, 81 R50

\section{Introduction}

The study of non-commutative geometry on quantum groups requires
to select a {\it first order differential calculus}
(that is, a FODC) on the corresponding coordinate Hopf algebra.
It is natural to assume that this differential calculus is compatible with
left and right translations of the quantum group. Such calculi are
called {\it bicovariant} (see Section \ref{ss-basics} for the precise
definition) and have been invented by S.\,L.~Woronowicz in his seminal
paper \cite{a-Woro2}. In the meantime there exists a well developed general
theory of such calculi (see, for instance, \cite{b-KS}, Chapter 14).

Whereas there exists a unique distinguished "classical"
differential calculus on any $C^\infty$-manifold, no functorial method
is known in order to get a canonical differential calculus on a
quantum group. Thus the problem of classifying all possible finite
dimensional bicovariant FODC on a Hopf algebra arises.
This problem was investigated under various additional assumptions
in \cite{a-SchSch1,a-SchSch2}, \cite{a-Majid2} and \cite{a-BGMST}.
The present paper gives
a complete solution of the classification problem for the coordinate
Hopf algebras $\OGq$ of the
quantum groups $G_q = SL_q(n+1),Sp_q(2n)$
when the parameter $q$ is transcendental. Our main result
(Theorem \ref{t-main}) says that up to forming direct sums
all finite dimensional bicovariant FODC on
$\mc{O}(SL_q(n+1))$ and $\mc{O}(Sp_q(2n))$ can
be obtained by the same general method. In special examples, this
method was first used in [4] and [10]. It was extended to arbitrary
coquasitriangular Hopf algebras in [11], 14.5-6, where the structure
of such calculi is developed in detail.
More precisely, any {\it irreducible}
finite dimensional bicovariant FODC is parametrized by a pair
$(v,\zeta)$ of an irreducible matrix corepresentation
$v=(v_{ij})_{i,j=1,\ldots,k}$ of $G_q=SL_q(n+1),Sp_q(2n)$ and a complex
number $\zeta$ such that $\zeta^{n+1}=1$ for $SL_q(n+1)$ and $\zeta^2=1$
for $Sp_q(2n)$. The dimension of the corresponding FODC is
$k^2$. A remarkable byproduct of the classification theorem
is that any finite dimensional bicovariant FODC is
inner, that is, there is a biinvariant one form $\theta$ such that
$\dif a=\theta a-a\theta$ for all $a\in\OGq$. Another
important corollary is that the quantum Lie algebra of any such calculus
can be obtained by means of a central element of the Hopf dual $\OGqn$
(see \ref{ss-centrgenFODC}).


The proof of our main classification theorem \ref{t-main} relies essentially
on two deep algebraic results. The first one is a
theorem due to A.~Joseph and G.~Letzter \cite{a-JoLet2} on
the locally finite part for the left adjoint action of the corresponding
quantized universal enveloping algebra $\Uqg$ (see \ref{ss-locfinuqg}) and
the second is a result of A.~Joseph \cite{b-Joseph} on the description of
the Hopf dual of $\RGq$ (see \ref{ss-hopfdual}).

This paper is organized as follows. Section \ref{sec-bicDC} contains
general notions and facts on bicovariant FODC which are needed later.
Proposition \ref{s-qliealg1} relates the bicovariant FODC on a Hopf algebra
$\A$ to $\adR$-invariant right coideals of the Hopf dual $\A^\circ$.
Moreover, the general method for the construction of such calculi from
\cite{b-KS} is briefly recalled. In Section \ref{sec-adrogq}
for an arbitrary simple Lie algebra the dual of the Hopf algebra $\RGq$
(see Proposition \ref{s-hopfdual})
and its finite dimensional
$\adR$-invariant coideals (see Theorem \ref{t-adinvkoid}) are
described when $q$ is transcendental. The latter result is of
interest in itself and it provides a classification of all
quantum Lie algebras of finite
dimensional bicovariant differential calculi on the Hopf algebra $\RGq$.
In Section \ref{sec-klassifikation} the main results of this paper
(Theorem \ref{t-main} and Corollaries
\ref{f-FODCinner}--\ref{f-FODCirred}) on the
classification problem for the Hopf algebras
$\mc{O}(SL_q(n+1))$ and $\mc{O}(Sp_q(2n))$ are stated and proved.
It is shown therein that all finite dimensional bicovariant
calculi can be given by the general methods
developed in 2.3 and 2.4. As
a byproduct we prove that the coquasitriangular
Hopf algebras $\mc{O}(SL_q(n+1))$ and $\mc{O}(Sp_q(2n))$ are factorizable.

In the formulation given in \ref{ss-mainresults} the classification results
are no longer valid for the quantum groups $\mc{O}(O_q(N))$.
For instance, the two four dimensional
irreducible bicovariant FODC on the Hopf algebra $\mc{O}(O_q(3))$
constructed in \cite{a-SchSch2} are not of the form
$\Gamma_\zeta (v)$ for some corepresentation $v$ of $O_q(3)$. The reason
for this failure lies in the fact that the tensor product
$\mc{O}(O_q(N))\otimes\mathbb{CZ}^n_2$ is not the full Hopf dual of
$\mathbb{CZ}^2\smashright\mc{U}_q(\lieso_N)$, because it does not
contain the matrix coefficients
of the spinor representations.
However, large parts of the proofs given in this paper carry over
with some modifications to the orthogonal case. In fact, if the
coordinate Hopf algebra $\mc{O}(O_q(N))$ is replaced by the
larger Hopf algebra generated by the matrix coefficients of the spinor
representations, then the results in \ref{ss-mainresults} are also valid
for this Hopf algebra.

We close this introduction by fixing some assumptions and
notations that
are used in the sequel. All algebras and Lie algebras are
over the complex field.
We shall use the convention to sum over repeated indices belonging to
different terms. The symbol $\rule{1ex}{1ex}$ marks the end of a proof,
an example or a remark.
Let $\A$ be a Hopf algebra. The
comultiplication, the counit and the antipode of $\A$ are denoted by
$\Delta$, $\varepsilon$ and $S$, respectively. Throughout we shall use the
Sweedler notation in the form $\Delta(a)=a_{(1)}\otimes a_{(2)}$. By a
\textit{right coideal} of $\A$ we mean a linear subspace $\mc{C}$ of
$\A$ such that $\Delta(\mc{C})\subseteq\mc{C}\otimes\A$. As
usual, $\A^\circ$ denotes the Hopf dual of $\A$. For $a,b\in\A$, let
$\adR(a)b=S(a_{(1)})ba_{(2)}$. Then the map $a\to\adR(a)$ is a
representation of the algebra $\A$ on itself, called the \textit{right
adjoint action} of $\A$. The map $\AdR:\A\to\A\otimes\A$ defined by
$\AdR(a)=a_{(2)}\otimes S(a_{(1)})a_{(3)}$, $a\in\A$, is the right
adjoint coaction of $\A$. A matrix $v=(v_j^i)_{i,j=1,\ldots,m}$ of
elements $v^i_j\in\A$ is called a \textit{matrix corepresentation}
of $\A$ if
$\varepsilon(v^i_j)=\delta_{ij}$ and $\Delta(v^i_j)=v^i_k\otimes v^k_j$
for $i,j=1,\ldots,m$. For a corepresentation $v$ of $\A$,
let $\mc{C}(v)$ denote the subcoalgebra of $\A$ spanned by the
matrix elements of $v$.

The main results of this paper were presented at the
quantum group conference in Prague (June 1997)
and an earlier version appeared in $q$-alg/9707032.

\section{Construction of bicovariant first order differential calculi}
\label{sec-bicDC}

\subsection{Bicovariant first order differential caluli: basic concepts}
\label{ss-basics}

We shall use the general framework of bicovariant differential calculus
developed by S.\,L.~Woronowicz \cite{a-Woro2}, see also \cite{b-KS},
Chapter 14. In this subsection we collect the main notions and facts needed
in what follows. Suppose that $\A$ is a Hopf algebra.

A \textit{first order differential calculus} (abbreviated, a FODC) over
$\A$ is an $\A$-bimodule $\Gamma$ equipped with a linear mapping
$\dif:\A\rightarrow\Gamma$, called the differentiation of the FODC, such
that:

\begin{description}
\item[(i)] $\dif$ satisfies the Leibniz rule $\dif(ab)=a\cdot\dif b+
\dif a\cdot b$ for any $a,b\in\A$,
\item[(ii)] $\Gamma$ is the linear span of elements $a\cdot\dif b\cdot c$
with $a,b,c\in\A$.
\end{description}

A FODC $\Gamma$ over $\A$ is called \textit{bicovariant} if there exist
linear mappings $\Delta_L:\Gamma\rightarrow\A\otimes\Gamma$
and $\Delta_R:\Gamma\rightarrow\Gamma\otimes\A$ such that
\begin{equation}
\Delta_L(a\dif b)=\Delta(a)(\id\otimes\dif)\Delta(b)\quad\mbox{and}\quad
\Delta_R(a\dif b)=\Delta(a)(\dif\otimes\id)\Delta(b)\quad
\mbox{for all }
a,b\in\A.
\end{equation}

A bicovariant FODC $\Gamma$ is called \textit{inner} if there exists a
biinvariant one form $\theta\in\Gamma$ (that is,
$\Delta_L(\theta)=1\otimes\theta$ and $\Delta_R(\theta)=\theta\otimes 1$)
such that
\begin{equation} \label{eq-innerdif}
\dif a=\theta a-a\theta,\quad a\in\A.
\end{equation}

By the dimension of a bicovariant FODC we mean the dimension of the vector
space $\Gammalinv=\{\omega\in\Gamma\,|\,\Delta_L(\omega)=1\otimes\omega\}$
of left invariant one forms. All bicovariant FODC occurring in this paper
are assumed to be \textit{finite dimensional}.

Let $\Gamma$ be a bicovariant FODC over $\A$. Then the set
\begin{equation}
\mc{R}_\Gamma:=\{a\in\ker\varepsilon\subset\A\,|\,
S(a_{(1)})\dif a_{(2)}=0\}
\end{equation}
is an $\AdR$-invariant right ideal of $\ker\varepsilon$.
Conversely, for any $\AdR$-invariant right ideal $\mc{R}$ of
$\ker\varepsilon$ there exists a bicovariant FODC $\Gamma$ such that
$\mc{R}=\mc{R}_\Gamma$ (\cite{b-KS}, Proposition 14.7). The linear
subspace
\begin{equation}
\X_\Gamma:=\{X\in\A^\prime\,|\,X(1)=0\quad\mbox{and}\quad X(a)=0
\quad\mbox{for}\quad a\in\mc{R}_\Gamma\}
\end{equation}
is called the \textit{quantum Lie algebra} of $\Gamma$.

\begin{satz}\label{s-qliealg}
(i) $\X_\Gamma$ is an $\adR$-invariant vector space of the dual Hopf
algebra $\A^\circ$ satisfying $\Delta(X)-\varepsilon\otimes X\in
\X_\Gamma\otimes\A^\circ$ for all $X\in\X_\Gamma$.\\
(ii) If $\A^\circ$ separates the elements of $\A$, then any finite
dimensional $\adR$-invariant vector space $\X\subset\A^\circ$ satisfying
$\Delta(X)-\varepsilon\otimes X\in\X\otimes\A^\circ$ and $X(1)=0$
for all $X\in\X$ is the quantum Lie algebra of a unique bicovariant FODC
over $\A$.
\end{satz}

\begin{bew}{}
\cite{b-KS}, Corollary 14.10.
\end{bew}

For the considerations below the following characterization of quantum
Lie algebras of bicovariant FODC will be more convenient.
Let us define a projection $P_\varepsilon:\A^\circ\to\A^\circ$ by
$P_\varepsilon(f):=f-f(1)\varepsilon$. Obviously, $P_\varepsilon(\A^\circ)=
\{f\in\A^\circ\,|\,f(1)=0\}$. For a linear subspace $\X$ of $\A^\circ$ we
set $\bar{\X}:=\X+\mathbb{C}\varepsilon$.

\begin{satz} \label{s-qliealg1}
(i) If $\X_\Gamma$ is the quantum Lie algebra of a bicovariant FODC, then
$\bar{\X}_\Gamma$ is an $\adR$-invariant right coideal of $\A^\circ$.\\
(ii) Suppose that $\A^\circ$ separates the elements of $\A$.
If $\bar{\X}$ is a finite dimensional $\adR$-invariant right coideal of
$\A^\circ$ containing $\varepsilon$, then $\X:=P_\varepsilon(\bar{\X})$ is
the quantum Lie algebra of a unique bicovariant FODC over $\A$.
\end{satz}

\begin{bew}{}
(i) We apply Proposition \ref{s-qliealg}(i). Because $\X_\Gamma$ is
$\adR$-invariant, $\bar{\X}_\Gamma$ is $\adR$-invariant as well.
Since $\Delta(\varepsilon)=\varepsilon\otimes\varepsilon$
and $\Delta(X)\in\varepsilon\otimes X+\X_\Gamma\otimes\A^\circ\subset
\bar{\X}_\Gamma\otimes\A^\circ$ for $X\in\X_\Gamma$, we conclude that
$\bar{\X}_\Gamma$ is a right coideal of $\A^\circ$.\\
(ii) For $Y\in\bar{\X}$ and $f\in\A^\circ$, we have
$\adR(f)(Y-Y(1)\varepsilon)=\adR(f)Y-f(1)Y(1)\varepsilon\in\bar{\X}$ and
$\adR(f)(Y-Y(1)\varepsilon)(1)=0$. Hence $\X$ is $\adR$-invariant.
Further, for $X\in\X$ we compute
$\Delta(X)-\varepsilon\otimes X=
\Delta(X)-\varepsilon\otimes X_{(1)}(1)X_{(2)}=
(P_\varepsilon\otimes\id)(\Delta(X)-\varepsilon\otimes X)\in
\X\otimes\A^\circ.$
Since obviously $X(1)=0$ for $X\in\X$, the assertion follows from
Proposition \ref{s-qliealg}(ii).
\end{bew}

Suppose that $\mc{N}$ is a subbimodule of the $\A$-bimodule $\Gamma$ such
that $\Delta_L(\mc{N})\subseteq\A\otimes\mc{N}$ and
$\Delta_R(\mc{N})\subseteq\mc{N}\otimes\A$.
Let $\pi:\Gamma\to\Gamma/\mc{N}$ denote the canonical map. Then the quotient
$\A$-bimodule $\Gamma/\mc{N}$ endowed with the linear mapping
$\tilde{\dif}:=\pi\circ\dif$ is also a bicovariant FODC over $\A$, called
a bicovariant quotient FODC of $\Gamma$. We shall say that the bicovariant
FODC $\Gamma$ is \textit{irreducible} if
$\Gamma$ has no nontrivial bicovariant quotient FODC, that is, if there is
no $\A$-subbimodule $\mc{N}\neq\{0\},\Gamma$ of $\Gamma$ such that
$\Delta_L(\mc{N})\subseteq\A\otimes\mc{N}$ and
$\Delta_R(\mc{N})\subseteq\mc{N}\otimes\A$.
An $\adR$-invariant right coideal $\X$ of a Hopf algebra is called
\textit{irreducible} if it does not contain
an $\adR$-invariant right coideal $\Y$ such that $\Y\not=\{0\},\X$.
One easily verifies the following lemma.

\begin{lemma}\label{l-equivirredFODC}
For a bicovariant FODC $\Gamma$ the following
statements are equivalent:\\
(i) The FODC $\Gamma$ is irreducible.\\
(ii) There is no $\adR$-invariant right subcoideal $\mc{Y}\subset
\bar{\X}_\Gamma$ of $\A^\circ$ containing $\varepsilon$ such that
$\mc{Y}\not=\mathbb{C}\varepsilon,\bar{\X}_\Gamma$.\\
(iii) There is no $\AdR$-invariant right ideal $\mc{I}$ of $\A$ such
that $\mc{R}_\Gamma\subset\mc{I}\subset\ker\varepsilon$,
$\mc{I}\not=\mc{R}_\Gamma$ and $\mc{I}\not=\ker\varepsilon$.
\end{lemma}

For $i\in\{1,\ldots,k\}$, let $\Gamma_i$ be a FODC over
$\A$ with differentiation $\dif_i$. Define a linear mapping $\dif=
\dif_1+\cdots+\dif_k$ of $\A$ into the direct sum $\tilde\Gamma:=
\Gamma_1\oplus\cdots\oplus\Gamma_k$ of $\A$-bimodules
$\Gamma_1,\ldots,\Gamma_k$. Clearly, $\Gamma:=\A{\cdot}\dif\A{\cdot}\A$
is a FODC over $\A$ with differentiation $\dif$.
We call this FODC $\Gamma$ the \textit{sum} of the FODC
$\Gamma_1,\ldots,\Gamma_k$. Suppose that the FODC
$\Gamma_1,\ldots,\Gamma_k$ are bicovariant. Then it is not difficult
to show that $\Gamma$ is also bicovariant and that
$\X_\Gamma=\X_{\Gamma_1}+\cdots+\X_{\Gamma_k}$ and
$\mc{R}_{\Gamma}=\mc{R}_{\Gamma_1}\cap\ldots\cap\mc{R}_{\Gamma_k}$. We
shall say that $\Gamma$ is a \textit{direct sum} of the FODC
$\Gamma_1,\ldots,\Gamma_k$ if $\tilde\Gamma=\Gamma$.
One easily verifies

\begin{lemma}\label{l-direktesumme}
$\Gamma$ is the direct sum of $\Gamma_1,\ldots,\Gamma_k$ if and only
if $\X_\Gamma$ is the direct sum of quantum Lie algebras
$\X_{\Gamma_1},\ldots,\X_{\Gamma_k}$ (that is, if
$\X_{\Gamma_i}\cap\X_{\Gamma_j}=\{0\}$ for $i\neq j$,
$i,j=1,\ldots,k$.)
\end{lemma}

Let $\Gamma_1$ and $\Gamma_2$ be two bicovariant FODC over $\A$. Suppose
that $\Gamma_i$ is a subbimodule of a possibly larger $\A$-bimodule
$\tilde{\Gamma}_i$ which contains a biinvariant element
$\theta_i\in\Gamma_i$ such that $\dif_ia=\theta_ia-a\theta_i,$ $a\in\A$
for $i=1,2$. Let $\dif$ be the linear mapping of $\A$ into the tensor
product $\tilde{\Gamma}=\tilde{\Gamma}_1\otimes_\A\tilde{\Gamma}_2$
of the two $\A$-bimodules $\tilde{\Gamma}_1$ and $\tilde{\Gamma}_2$
defined by
\begin{equation}\label{eq-innerdiftensor}
\dif a=\theta a-a\theta\equiv
\dif_1a\otimes_\A\theta_2+\theta_1\otimes_\A\dif_2a,\quad a\in\A,
\end{equation}
where $\theta:=\theta_1\otimes_\A\theta_2$. It is straightforward to
check that $\Gamma:=\A\cdot\dif\A\cdot\A$ is again a bicovariant FODC over
$\A$ with differentiation $\dif$. This FODC $\Gamma$ is called the
\textit{tensor product} of the FODC $\Gamma_1$ and $\Gamma_2$.

\subsection{Construction of bicovariant FODC on coquasitriangular Hopf
algebras} \label{ss-coquasiFODC}

In this subsection we briefly repeat the general method for the construction
of bicovariant FODC (see \cite{b-KS}, Sect.\ 14.5, for a detailled treatment).
Let $\A$ be a coquasitriangular Hopf algebra (see \cite{b-LT}, \cite{b-KS}
or \cite{b-Majid1}). That is,
$\A$ is Hopf algebra equipped with a linear functional $\br$ on $\A \otimes \A$ which
is invertible with respect to the convolution multiplication and satisfies the
following  conditions
\begin{gather*}
\br( ab \otimes c)= \br( a \otimes c_{(1)})\br( b\otimes c_{(2)}), \quad
\br( a\otimes bc)= \br( a_{(1)}\otimes c)\br( a_{(2)}\otimes b), \\
\br( a_{(1)}\otimes b_{(1)}) a_{(2)} b_{(2)}=\br( a_{(2)}\otimes
b_{(2)}) b_{(1)} a_{(1)}
\end{gather*}
for $a,b,c \in \A$. Such a linear form $\br$ is called a universal $r$-form of
the Hopf algebra $\A$.

Suppose that $\br$ is a fixed universal $r$-form  and
$v=(v^i_j)_{i,j=1,\ldots,m}$ is a matrix corepresentation of
$\A$. Let $\bar{\br}$ denote the convolution inverse of $\br$. We define linear
functionals ${l^{\pm}}^i_j,l^i_j\in\A^\circ$ by
\begin{equation}\label{eq-generlfunc}
{l^+}^i_j(\cdot)=\br(\cdot\otimes v^i_j),\qquad
{l^-}^i_j(\cdot)=\bar{\br}(v^i_j\otimes\cdot),\qquad
l^i_j=S({l^-}^i_k){l^+}^k_j.
\end{equation}

The $m\times m$-matrices $L^+$ and $L^{-,\contr}$ with entries
$(L^+)^i_j={l^+}^i_j$ and $(L^{-,\contr})^i_j=S({l^-}^j_i)$ are then
representations of the algebra $\A$ and the complex vector space
\begin{equation}\label{eq-lfunkraum}
\X\cont(v):=\mathrm{Lin}\,\{l^i_j\,|\,i,j=1,\ldots,m\}
\end{equation}
of the functionals $l^i_j$ is an $\adR$-invariant right coideal of
$\A^\circ$. Let $v\cont$ be the contragredient
matrix corepresentation of $\A$, that is, $v\cont$ is the
$m\times m$-matrix with entries $(v\cont)^i_j=S(v^j_i)$. Then the pairs
$\Gamma=(v,L^{-,\contr})$ and $\Gamma\cont=(v\cont,L^+)$ are bicovariant
bimodules over $\A$ (see \cite{b-KS}, Sect.\ 14.5, for details).
Hence the tensor product
\begin{equation}\label{eq-basicbikbim}
\tilde{\Gamma}:=\Gamma\otimes_\A\Gamma\cont=
(v\otimes v\cont,L^{-,\contr}\otimes L^+)
\end{equation}
of the bicovariant bimodules $\Gamma$ and $\Gamma\cont$ is again
a bicovariant bimodule. Let us fix $\{\theta_{ij}\}_{i,j=1,\ldots,m}$ be a
basis of the vector space $\linv{\tilde{\Gamma}}$ such that
$\Delta_R(\theta_{ij})=\theta_{kt}\otimes v^k_i(v\cont)^t_j$. Obviously,
$\theta:=\sum_i\theta_{ii}$ is a biinvariant element of $\tilde{\Gamma}$.
Equipped with the differentiation defined by (\ref{eq-innerdif}),
$\Gamma(v):=\A\cdot\dif\A\cdot\A$ becomes a bicovariant FODC $\Gamma(v)$
with quantum Lie algebra $\X(v)$ spanned by the linear functionals
\begin{equation}
X_{ij}=l^i_j-\delta_{ij}\varepsilon,\quad i,j=1,\ldots,m.
\end{equation}
If the functionals $X_{ij}$, $i,j=1,\ldots,m$, are linearly
independent, then we clearly have $\Gamma(v)=\tilde{\Gamma}$.

\begin{bsp}\label{b-trivFODC}
If $v=1$, then $\tilde{\Gamma}=(1,\varepsilon)$ and hence
$\X(v)=\{0\}$, so $\Gamma(v)\not=\tilde{\Gamma}$ and $\Gamma(v)$ is
the trivial bicovariant FODC with $\dif a=0$ for all $a\in\A$.
\end{bsp}

\subsection{Construction of bicovariant first order differential calculi
over the Hopf algebras $\mc{O}(SL_q(n+1))$ and $\mc{O}(Sp_q(2n))$}
\label{ss-OGqFODC}

Let us say  that a complex number $\zeta$ is \textit{admissible} if
$\zeta^{n+1}=1$ for $G_q=SL_q(n+1)$ and $\zeta^2=1$ for $G_q=Sp_q(2n)$.
For any admissible number $\zeta$ let $\varepsilon_\zeta$ denote the
multiplicative
linear functional on the algebra $\OGq$ such that
$\varepsilon_\zeta(u^i_j)=\zeta\delta^i_j$, $i,j=1,\ldots,N$.

First we define some 1-dimensional bicovariant FODC over $\OGq$.
For admissible $\zeta$, let $\Gamma_\zeta=(1,\varepsilon_\zeta)$ be the
one dimensional bicovariant bimodule over $\OGq$. That is, $\Gamma_\zeta$
has a free left module basis consisting of a single biinvariant element
$\theta_0$ and we have $\theta_0a=a_{(1)}\varepsilon_\zeta(a_{(2)})\theta_0$
for $a\in\OGq$. If $\zeta\not=1$, then $\Gamma_\zeta$ is an inner
bicovariant FODC over $\OGq$ with differentiation
(\ref{eq-innerdif}). It can be shown (see Remark 2.4 in \cite{a-SchSch1})
that any one dimensional bicovariant FODC over
$\OGq$ is of this form.

Now we specialize the procedure from the preceding section
to the coordinate Hopf algebra $\OGq$, where $G_q=SL_q(n+1)$ or
$G_q=Sp_q(2n)$. Recall (see \cite{b-KS}, 10.1.2) that these Hopf algebras
are coquasitriangular with universal $r$-form $\br$ such that
\begin{equation}\label{eq-rform}
\br(u^i_j\otimes u^n_m)=zR^{in}_{jm},\quad i,j,n,m=1,\ldots,N.
\end{equation}
Here $u^i_j$ are the matrix entries of the fundamental corepresentation of
$\OGq$, $R$ is the corresponding $R$-matrix (see \cite{a-FadResTak1},
(1.5) and (1.9), or \cite{b-KS}, (9.13) and (9.30)) and $z$ is a fixed
complex number such that $z^N=q^{-1}$ for $SL_q(N)$ and $z^2=1$ for
$Sp_q(N)$.

Suppose that $v=(v^i_j)_{i,j=1,\ldots,m}$ is a matrix corepresentation of
$\A=\OGq$ and $\zeta$ is an admissible number. Consider the bicovariant
bimodule
\begin{equation}\label{eq-tildeGamma}
\tilde{\Gamma}:=\Gamma_\zeta\otimes_\A\Gamma_v\otimes_\A\Gamma\cont_v=
(1\otimes v\otimes v\cont,\varepsilon_\zeta\otimes L^{-,\contr}\otimes L^+).
\end{equation}
The element $\theta:=\sum_i\theta_0\otimes_\A\theta_{ii}\in\tilde{\Gamma}$
is biinvariant and $\Gamma_\zeta(v):=\A\cdot\dif\A\cdot\A$ is a bicovariant
FODC over $\A=\OGq$ with differentiation (\ref{eq-innerdif}).
If $\zeta=1$, then $\Gamma_\zeta(v)$ is just the FODC $\Gamma(v)$ defined
in \ref{ss-coquasiFODC}. If $\zeta\not=1$, then the FODC $\Gamma_\zeta(v)$
is the tensor product of the one dimensional FODC $\Gamma_\zeta$ and the
FODC $\Gamma(v)$.

The quantum Lie algebra $\X_\zeta(v)$ of the bicovariant FODC
$\Gamma_\zeta(v)$ is then the linear span of the functionals
\begin{equation}\label{eq-ogqliealgfunc}
X_{ij}=\varepsilon_\zeta l^i_j-\delta_{ij}\varepsilon,\quad i,j=1,\ldots,m.
\end{equation}

\begin{satz}\label{s-tensorFODC}
Let $\zeta,\zeta'$ be admissible numbers for $G_q$ and
$v,w$ corepresentations of $\OGq$. Then the bicovariant FODC
$\Gamma_{\zeta\zeta'}(v\otimes w)$ and $\Gamma_\zeta(v)\otimes
\Gamma_{\zeta'}(w)$ (see \ref{ss-basics}) are isomorphic.
\end{satz}

\begin{bew}{}
Because of (\ref{eq-lfunkraum}) the sets
$\{l(ab)\,|\,a\in\mc{C}(v),b\in\mc{C}(w)\}$ and
$\{l(a)l(b)\,|\,a\in\mc{C}(v),b\in\mc{C}(w)\}$ span the
$\adR$-invariant right coideals $\X\cont(v\otimes w)$
and $\X\cont(v)\cdot\X\cont(w)$, respectively.
We compute
\begin{align*}
l(ab)&=
l(a_{(1)})\adR(l^+(a_{(2)}))l(b).
\end{align*}
Since $\Delta:\mc{C}(v)\to\mc{C}(v)\otimes \mc{C}(v)$ and $\X\cont(w)$
is $\adR$-invariant, we obtain that $l(ab)\in\X\cont(v)\cdot\X\cont(w)$.
On the other hand, by formula (10.29) in \cite{b-KS}, we have
\begin{displaymath}
l(a)l(b)=\br(b_{(1)}\otimes S(a_{(3)}))\br(b_{(3)}\otimes a_{(2)})l(a_{(1)}b_{(2)}).
\end{displaymath}
Hence it follows that $l(a)l(b)\in\X\cont(v\otimes w)$ for all
$a\in\mc{C}(v)$ and $b\in\mc{C}(w)$. This proves that
$\X\cont(v\otimes w)=\X\cont(v)\cdot\X\cont(w)$.

The quantum Lie algebra of the bicovariant FODC
$\Gamma_{\zeta\zeta'}(v\otimes w)$ is given by
(\ref{eq-ogqliealgfunc}). Clearly, we have
$\X_{\zeta\zeta'}(v\otimes w)=
P_\varepsilon(\varepsilon_{\zeta\zeta'}\X\cont(v\otimes w))$.
Let $\theta_1$ and $\theta_2$ be nonzero biinvariant elements in
$\Gamma_\zeta(v)$ and $\Gamma_{\zeta'}(w)$, respectively.
If one of the pairs $(v,\zeta)$ or $(w,\zeta')$ is equal to $(1,1)$,
then for this pair we take the one dimensional
bicovariant bimodule $\tilde{\Gamma}=(1,\varepsilon)$ given by
(\ref{eq-basicbikbim}) and a biinvariant element therein.
By the definition of the tensor product of FODC,
the biinvariant element $\theta_1\otimes_\A\theta_2$ defines the
differentiation of $\Gamma_\zeta(v)\otimes \Gamma_{\zeta'}(w)$
(see (\ref{eq-innerdiftensor})). Obviously, the quantum Lie
algebra of $\Gamma_\zeta(v)\otimes \Gamma_{\zeta'}(w)$ is
$P_\varepsilon(\varepsilon_{\zeta}\X\cont(v)
\varepsilon_{\zeta'}\X\cont(w))$. Since $\varepsilon_\zeta$ and
$\varepsilon_{\zeta'}$ are central
and $\X\cont(v\otimes w)=\X\cont(v)\cdot\X\cont(w)$ as shown in the preceding
paragraph, we have $P_\varepsilon(\varepsilon_{\zeta}\X\cont(v)
\varepsilon_{\zeta'}\X\cont(w))=
P_\varepsilon(\varepsilon_{\zeta\zeta'}\X\cont(v\otimes w))$.
That is, the quantum Lie algebras of the bicovariant FODC $\Gamma_{\zeta\zeta'}(v\otimes w)$ and
$\Gamma_{\zeta}(v)\otimes\Gamma_{\zeta'}(w)$ coincide.
Hence the FODC are isomorphic.
\end{bew}

Set $(D^{-1})^j_i:=\br(S^2(v^j_n)\otimes v^n_i)$ for $i,j=1,\ldots,m$.
Then the element
\begin{equation}\label{eq-zentrinGamma}
c_\zeta(v):=\varepsilon_\zeta\mathrm{Tr}\,LD^{-1}=
\sum\nolimits_{i,j}\varepsilon_\zeta l^i_j(D^{-1})^j_i.
\end{equation}
belongs to the centre of the Hopf dual $\A^\circ$
(\cite{b-KS}, Proposition 10.15(ii)), so each
quantum Lie algebra $\X_\zeta(v)$ contains the central element
\begin{displaymath}
P_\varepsilon(c_\zeta(v))=c_\zeta(v)-\sum\nolimits_i(D^{-1})^i_i\varepsilon .
\end{displaymath}

\subsection{Bicovariant FODC generated by central elements of $\A^\circ$}
\label{ss-centrgenFODC}

In this brief subsection we recall another general method for the
construction of bicovariant FODC which was invented in \cite{a-BrzMaj2},
see also \cite{a-Majid2}. We
state the corresponding result as

\begin{satz}\cite{a-BrzMaj2,a-Majid2}\label{s-centrgenFODC}
Suppose that $\A$ is a Hopf algebra such that its Hopf dual $\A^\circ$
separates the elements of $\A$. For any central element $c\in\A^\circ$
there exists a bicovariant FODC $\Gamma(c)$ over $\A$ with quantum Lie
algebra
\begin{equation}\label{eq-zentrgener}
\X[c]=\mathrm{Lin}\,\{\chi_a:=c_{(2)}(a)c_{(1)}-c(a)\varepsilon\,|\,
a\in\A\}.
\end{equation}
\end{satz}
\begin{bew}{}
See \cite{b-KS}, Proposition 14.11.
\end{bew}

Note that (\ref{eq-zentrgener}) differs from the corresponding
formula in \cite{a-Majid2}, because we use another
definition of the quantum Lie algebra.

If we write $\Delta(c)=\sum_{i}x_i\otimes y_i$ with $\{y_i\}$ linearly
independent, then $\X[c]$ is obviously equal to the linear span of
functionals $x_i-x_i(1)\varepsilon$. In particular, we see that the
vector space $\X[c]$ and hence the bicovariant FODC $\Gamma(c)$
are finite dimensional.

\section{The right adjoint action of the Hopf dual $\RGqn$}
\label{sec-adrogq}

In this section $\mf{g}$ denotes a simple Lie algebra with
Cartan matrix $(a_{ij})$ and rank $n$ and $(d_i\delta_{ij})$ is a
diagonal matrix such that
$(d_ia_{ij})$ is symmetric. Throughout we
suppose that $q$ is a transcendental complex number.
(However, some assertions and
arguments require only that $q$ is
not root of unity.)

\subsection{The extended Drinfeld-Jimbo algebra $\Uqg$}
\label{ss-Uqg}

Let $\Uqg$ denote the unital complex algebra with
generators $E_i$, $F_i$, $K_i$, $K_i^{-1}$, $i=1,\ldots,n$,
subject to the relations
\begin{gather}
K_iK_j=K_jK_i,\quad K_iK_i^{-1}=K_i^{-1}K_i=1,\notag\\
\label{eq-kommutEK}
K_iE_jK_i^{-1}=q_j^{\delta^i_j}E_j,\quad
K_iF_jK_i^{-1}=q_j^{-\delta^i_j}F_j,\notag\\
\label{eq-kommutEF}
E_iF_j-F_jE_i=
\delta^i_j\frac{K_{\alpha_i}-K_{\alpha_i}^{-1}}{q_i-q_i^{-1}},\\
\sum_{r=0}^{1-a_{ij}}(-1)^r\left[{1-a_{ij}\atop r}\right]_{q_i}
E_i^{1-a_{ij}-r}E_jE_i^r=0,\quad i\not=j,\notag\\
\sum_{r=0}^{1-a_{ij}}(-1)^r\left[{1-a_{ij}\atop r}\right]_{q_i}
F_i^{1-a_{ij}-r}F_jF_i^r=0,\quad i\not=j,\notag
\end{gather}
where
\begin{displaymath}
\left[{m\atop k}\right]_q=\frac{[m]_q!}{[k]_q![m-k]_q!},\quad
[m]_q!=[m]_q[m-1]_q\cdots[1]_q,\quad [0]_q!=1\mbox{ and }
[m]_q=\frac{q^m-q^{-m}}{q-q^{-1}}.
\end{displaymath}
In equation (\ref{eq-kommutEF}) the symbol $K_{\alpha_i}$ denotes the
element $K_{\alpha_i}=\prod_{j=1}^nK_j^{a_{ji}}$ and
$K_{\alpha_i}^{-1}$ its inverse. The generators $K_i$ correspond to the
fundamental weights
$\fuw_i$ of the Lie algebra $\mf{g}$. As usual, $\Uqnp$ and $\Uqnm$ are
the subalgebras of $\Uqg$ generated by the elements
$E_i$ resp. $F_i$, $i=1,\dots ,n$.

There is a Hopf algebra structure on $\Uqg$ determined by the formulas
\begin{gather*}
\Delta(K_i)=K_i\otimes K_i,\quad
\Delta(K_i^{-1})=K_i^{-1}\otimes K_i^{-1},\\
\Delta(E_i)=E_i\otimes K_{\alpha_i}+1\otimes E_i,\quad
\Delta(F_i)=F_i\otimes 1+K_{\alpha_i}^{-1}\otimes F_i,\\
\varepsilon(K_i)=\varepsilon(K_i^{-1})=1,\quad
\varepsilon(E_i)=\varepsilon(F_i)=0,\\
S(K_i)=K_i^{-1},\quad S(K_i^{-1})=K_i,\quad
S(E_i)=-E_iK_{\alpha_i}^{-1},\quad S(F_i)=-K_{\alpha_i}F_i.
\end{gather*}

The Hopf subalgebra $\Uqg_{DJ}$ of $\Uqg$
generated by the elements $E_i$, $F_i$ and $K_{\alpha_i}$, $i=1,\ldots,n,$
is isomorphic to the "usual" Drinfeld-Jimbo algebra.
For the considerations in Section 4 it is crucial to work with the extended
Drinfeld-Jimbo algebra $\Uqg$, because this algebra is also generated by
the so-called $L$-functionals ${l^\pm}^i_j$ (see
\cite{a-FadResTak1} or \cite{b-KS}, 8.5.3, for details).


\subsection{The dual Hopf algebra $\RGqn$}
\label{ss-hopfdual}



Let $\RGq$ denote the Hopf subalgebra of $\Uqg^\circ$ spanned by the
matrix coefficients of all
finite dimensional type 1 highest weight representations of $\Uqg$. The
purpose of this subsection is to repeat a result of Joseph [7] on the
description of the Hopf dual $\RGqn$.

Let $\chru$ be the multiplicative group in $\Uqg$ generated by
$K_i$, $i=1,\ldots,n$. First we describe the Hopf dual of the group algebra
$\mathbb{C}\chru$ (see \cite{b-Joseph}, Lemma 9.4.8(i)). Let
$(\mathbb{C}^\times)^n$ be the product of $n$
copies of the multiplicative group of nonzero complex numbers. For any
$\mu=(\mu_1,\ldots,\mu_n)\in(\mathbb{C}^\times)^n$, there is a
unique multiplicative linear functional
$f_\mu$ on $\mathbb{C}\chru$ such that $f_\mu(K_i)=\mu_i$,
$i=1,\ldots,n$. Hence there are numbers $\tilde{\mu}_i\neq 0$ such that
$f_\mu(K_{\alpha_i}^{-1})=\tilde{\mu}_i$.
Further, for any $i=1,\ldots,n$ there is a unique linear functional
$g_i$ on $\mathbb{C}\chru$ such that $g_i(K_j)=(a^{-1})_{ij}$ and
$g_i(ab)=\varepsilon(a)g_i(b)+g_i(a)\varepsilon(b)$
for all $a,b\in\mathbb{C}\chru$. Here $(a^{-1})_{ij}$ are the
entries of the inverse Cartan matrix. From the
normalization $g_i(K_j)=(a^{-1})_{ij}$ we get
$g_i(K_{\alpha_j})=\delta^i_j$ which implies the
simple commutation rules (\ref{eq-kommutfg}) below
between $E_i, F_i$ and $g_j$.
Now the Hopf dual
$\U_0:=(\mathbb{C}\chru)^\circ$ is the (commutative
and cocommutative) Hopf algebra generated by the functionals
$f_\mu$ and $g_i$ with relations
$f_\mu f_{\mu'}=f_{\mu\mu'}, g_ig_j=g_jg_i, f_\mu g_i=g_if_\mu$
and Hopf algebra structure given by
\begin{xalignat*}{3}
\Delta(f_\mu)&=f_\mu\otimes f_\mu,&
\varepsilon(f_\mu)&=1,&
S(f_\mu)&=f_{\mu^{-1}},\\
\Delta(g_i)&=g_i\otimes 1+1\otimes g_i,&
\varepsilon(g_i)&=0,&
S(g_i)&=-g_i.
\end{xalignat*}
Here $\mu^{-1}$ denotes the $n$-tuple $(\mu_1^{-1},\ldots,\mu_n^{-1})
\in(\mathbb{C}^\times)^n$ if $\mu=(\mu_1,\ldots,\mu_n)$.

As shown in \cite{b-Joseph}, the Hopf algebra $\U_0$ can be
embedded into the Hopf algebra $\RGqn$. For notational simplicity let
us identify $\U_0$ with its image in $\RGqn$. We
shall not need the explicit form of this embedding. It suffices to
know the following commutation relations
(see \cite{b-Joseph}, 9.4.8 and 9.4.9) of the functionals
$f_\mu, g_i \in \U_0$ with the
generators $E_i, F_i \in \Uqg$ in the algebra $\RGqn$:
\begin{equation}\label{eq-kommutfg}
f_\mu E_i={\tilde \mu_i}E_i f_\mu ,\quad
F_i f_\mu ={\tilde \mu_i} f_\mu F_i,\quad
E_ig_j=(g_j+\delta_j^i)E_i,\quad
g_jF_i=F_i(g_j+\delta_j^i).
\end{equation}

\begin{satz}\label{s-hopfdual}\cite{b-Joseph}
The multiplication map \quad
$\Uqnp \otimes \U_0\otimes \Uqnm \to \RGqn$
is an isomorphism of vector spaces.
\end{satz}

\begin{bew}{} Proposition 9.4.9 in \cite{b-Joseph}.
\end{bew}



\subsection{The locally finite part for the right adjoint action of
$\Uqg$}
\label{ss-locfinuqg}

In this subsection we restate the fundamental result of Joseph
and Letzter \cite{a-JoLet1,a-JoLet2,b-Joseph} about the locally finite
part of $\Uqg$. If $\B$ is a subalgebra of a Hopf algebra $\A$, we define
\begin{displaymath}
\mc{F}_\B(\A):=\{a\in\A\,|\,\mathrm{dim}\,\adR(\B)a<\infty\}.
\end{displaymath}
The vector space $\mc{F}(\A):=\mc{F}_\A(\A)$ is called
the \textit{locally finite part} of the Hopf algebra $\A$. We denote by
$P(\mf{g})=\{\lambda\in\mf{h}^*\,|\,\langle\alpha_i^\vee,
\lambda\rangle\in\mathbb{Z}\}$
the weight lattice of $\mf{g}$ and by
$P_+(\mf{g})$ the subset of dominant integral weights.
Let $\fuw_i\in P_+(\mf{g})$
be given by $\langle\alpha_i^\vee,\fuw_j\rangle=\delta_{ij}$.
There is a unique group
homomorphism $\tau:P(\mf{g})\to \chru$ such that
$\tau(\fuw_i)=K_i$. Further, let $\pchru$ be the
subsemigroup of $\chru$
generated by the elements $K_i^{2}=\tau(2\fuw_i)$, $i=1,\dots,n$, and $\tUqg$ the subalgebra of $\RGqn$ generated
by the elements $E_i, F_i$ and $f_\mu$.

\begin{satz}
\label{s-locfinpart}
(i) $\mc{F}(\Uqg)=\bigoplus_{\lambda\in P_+(\mf{g})}
\adR(\Uqg)\tau(-2\lambda)$. Further, the subspace
$W(\lambda):=\adR(\Uqg)\tau(-2\lambda)$ has dimension
$(\mathrm dim\,V(\lambda))^2$, where $V(\lambda)$
is the irreducible representation of $\Uqg$ with
highest weight $\lambda$ with respect to the
simple roots $\{-\alpha_1,\ldots,-\alpha_n\}$.

(ii) There is a subset $M$ of $(\comp^\times)^n$ such that
$\tUqg = \bigoplus_{\mu \in M} \mc{F}(\Uqg)\pchru f_\mu $.
\end{satz}

\begin{bew}{}
Since a subspace $E$ of $\Uqg$ is $\adL$-invariant if and
only if $S^{-1}(E)$ is $\adR$-invariant, (i) is only a restatement of
the corresponding result for the left adjoint action $\adL$ obtained in
\cite{a-JoLet2}. (ii) follows from \cite{b-Joseph}, 7.1.13.
\end{bew}

\subsection{Finite dimensional $\adR$-invariant right coideals of
$\RGqn$}

Let $C$ denote the set of all multiplicative linear functionals in the centre
of $\RGqn$.
Recall we assumed in this section
that $q$ is a transcendental complex number.

\begin{lemma}\label{l-adREik}
Suppose that $u\in \mc{F}(\Uqg)$, $\lambda \in P_+(\mf{g})$ and
$\mu\in(\comp^\times)^n$. If there is a natural number $k$ such that
\begin{equation}\label{eq-annu}
\adR(E_i^k)(u\tau(2\lambda)f_\mu) = 0
\end{equation}
for all i, then it follows that $u=0$ or
$\tau(2\lambda)f_\mu=c\tau(-2\gamma)$ for some $c\in C$ and
$\gamma \in P(\mf{g})$.
\end{lemma}

\begin{bew}{}
Since $u\in \mc{F}(\Uqg)$, there exists a natural number $l$ such that
$\adR(E_i^l)u = 0$ for all i. Let us fix the index $i$ and abbreviate
\begin{equation*}
b_{rs}:=q_i^{-s(l+k+r-s)}\left[{l+k+r\atop s}\right]_{q_i}
\prod_{j=0}^{l+k+r-s-1}
(\tilde{\mu}_iq_i^{2\langle \alpha_i^\vee,\lambda\rangle }-q_i^{-2j}).
\end{equation*}
{}From the relations $\adR(E_i^{l+k+r})(u\tau(2\lambda)f_\mu) = 0$
by (\ref{eq-annu}) and $\adR(E_i^l)u = 0$ we obtain
\begin{equation}\label{eq-sum}
\sum_{s=0}^{l-1}b_{rs}\adR(E_i^s)u \cdot E_i^{l+k+r-s}
\tau(2\lambda)f_\mu=0,\, r=0,\dots,l-1.
\end{equation}
By induction one proves that the determinant of the
matrix $(b_{rs})_{r,s=0,\dots,l-1}$ is equal to
\begin{equation}\label{eq-det}
\mathrm{det}\, (b_{rs})=q_i^{-l(l-1)(k+l)}
\prod_{j=1}^{l-1}
(\tilde{\mu}_iq_i^{2\langle \alpha_i^\vee,\lambda\rangle }-q_i^{2(l-j)})^j
(\tilde{\mu}_iq_i^{2\langle \alpha_i^\vee,\lambda\rangle }-q_i^{-2(j+k)})^{l-j}
\prod_{j=0}^k
(\tilde{\mu}_iq_i^{2\langle \alpha_i^\vee,\lambda\rangle }-q_i^{-2j})^l.
\end{equation}
Therefore, if $\mathrm{det}\, (b_{rs})\neq 0$, then it follows from
(\ref{eq-sum})
that $\adR(E_i^s)u \cdot E_i^{l+k+r-s}\tau(2\lambda)f_\mu=0$,
$r,s=0,\dots,l-1$.
In the case $s=0$ we get $u \cdot E_i^{l+k+r}\tau(2\lambda)f_\mu=0$
and hence $u=0$. If $\mathrm{det}\, (b_{rs})=0$, then
(\ref{eq-det})
implies that there are integers $n_i$ such that
$\tilde{\mu}_iq_i^{2\langle \alpha_i^\vee,\lambda\rangle }=q_i^{-2n_i}$
for all i. Set $\gamma:=\sum_i n_i\omega_i$.
From the commutation rules in the algebra $\tUqg$
we see that $c:=f_\mu \tau(2\lambda)\tau(2\gamma)$ commutes with all
generators
$E_i, F_i, f_\nu$, so $c$ is in the centre of $\tUqg$. Obviously,
$c$ is a character, because $f_\mu$, $\tau(2\lambda)$ and $\tau(2\gamma)$
are so. Hence we have $c\in C$ and $\tau(2\lambda)f_\mu=c\tau(-2\gamma)$.
\end{bew}

An immediate conseqence of Lemma \ref{l-adREik} is

\begin{folg}\label{f-adREiinj}
Suppose that $\mu\in(\comp^\times)^n$. If there is no $c\in C$ such that
$\mc{F}(\Uqg)\pchru f_\mu \subseteq c\mc{F}(\Uqg)\pchru$, then $\adR(E_i)$
is injective on the $\adR$-invariant subspace $\mc{F}(\Uqg)\pchru f_\mu$ for all $i$.
\end{folg}

Our next aim is to prove Proposition \ref{s-FRGq} below
which was suggested by the referee. The authors would like to thank the referee for his advice.

\begin{lemma}\label{l-FtUqg}
$\mc{F}(\tUqg)\subseteq\mc{F}(\Uqg)C$.
\end{lemma}

\begin{bew}{}
We first note that because $q$
is not a root of unity the $\adR$-right module $\mc{F}_{\Uqg}(\tUqg)$
of the algebra $\Uqg$ decomposes into a direct sum of simple submodules.
Let $V$ be such a submodule. Since the algebra $\tUqg$ is the sum of
$\adR$-invariant subspaces $\mc{F}(\Uqg)\pchru f_\mu$ by
Proposition \ref{s-locfinpart}, there is
$\mu\in(\comp^\times)^n$ such that
$V\subseteq \mc{F}(\Uqg)\pchru f_\mu$. Any element $v\in V$ is of
the form $v=u\tau(2\lambda)f_\mu $, where $u\in\mc{F}(\Uqg)$ and
$\lambda\in P_+(\mf{g})$, and satisfies the assumptions of Lemma
\ref{l-adREik}. If $u=0$, then $v=0$, so $v$ is trivially in $\mc{F}(\Uqg)C$.
If $u\not=0$, then $v=u\tau (2\lambda)f_\mu =uc\tau(-2\gamma) \in c\,\Uqg$
by Lemma \ref{l-adREik}. Since $v\in\mc{F}(\tUqg)$ by assumption, we get
$v \in \mc{F}(\Uqg)C$.
\end{bew}{}

Let $\mc{G}$ be the subalgebra of $\RGqn$ generated by the elements
$g_i$, $i=1,\ldots,n$. The elements
$\{g^\mf{k}=g_1^{k_1}\cdot\ldots\cdot g_n^{k_n}\,|\,
k_1,\ldots,k_n\in\mathbb{N}_0\}$, $\mf{k}=(k_1,\ldots,k_n)$,
form a vector space basis of $\mc{G}$. By Proposition \ref{s-hopfdual},
the multiplication map
gives a vector space isomorphism of $\tUqg\otimes\mc{G}$ and $\RGqn$.
That is, any $v\in \RGqn$ is a finite sum
$v=\sum_{\mf{k}} u_{\mf{k}} g^{\mf{k}}$ with uniquely determined elements
$u_{\mf{k}} \in \tUqg$. If $v\neq 0$, then the largest number
$|{\mf{k}}|:=k_1+\cdots+k_n$ such that
$u_{\mf{k}}\neq 0$ is called the degree of $v$ and denoted by $|v|$. We shall
write $v\simeq w$ if the degree of $v-w$ is lower
than the degree of $v$ and $w$.

\begin{lemma}\label{l-RGq}
Let $V$ be a simple submodule of $\RGqn$ with respect to the right adjoint
action $\adR$ of the algebra $\Uqg$.
Then we have $V\subset\tUqg$.
\end{lemma}

\begin{bew}{}
Assume the contrary. Then there exists an element $v\in V$ such that $|v|>0$.
Let us write $v$ in the form $v=\sum u_{\mf{k}}g^{\mf{k}} +
\sum u_{\mf{h}}g^{\mf{h}}$, where ${\mf{k}}$ and ${\mf{h}}$ run over all
multindices such that $|{\mf{k}}|=|v|$ and $|{\mf{h}}|<|v|$, respectively.
One easily verifies that
\begin{equation}\label{eq-sim}
\adR(f)(u_{\mf{k}}g^{\mf{k}})\simeq
(\adR(f)u_{\mf{k}})g^{\mf{k}}
\end{equation}
for all ${\mf{k}}$ and $f\in\Uqg$.
Hence $V_{\mf{k}}:=\adR(\Uqg)(u_{\mf{k}})$ is an
$\adR$-invariant finite dimensional vector space and
$u_{\mf{k}}\in \mc{F}(\tUqg)$ for all $\mf{k}$. Since
$\adR(E_i)$ raises the degree with respect to the usual $\mathbb{Z}$-gradation
(that is, $\mathrm{deg}E_i =-\mathrm{deg}F_i=1$ and
$\mathrm{deg}f_\mu =\mathrm{deg}g_i=0$) by one and the spaces $V_{\mf{k}}$
are finite dimensional, we can choose
$v\in V$ such that  $\adR(E_i)u_{\mf{k}}=0$ for all $\mf{k}$ and
all $i$. From the latter and (\ref{eq-sim}) we conclude that
$|\adR(E_i)v|<|v|$.
Therefore, $W_i:= \adR(\Uqg)(\adR(E_i)v)$ is a
submodule of the simple module $V$
such that $W'_i\neq V$. Thus, $W'_i= \{ 0 \}$ and hence $\adR(E_i)v=0$ for all $i$.

Fix a $\mf{k}'$ such that $u_{\mf{k}'}\neq 0$.
Since $|v|=|\mf{k}'|>0$, there exists
an index $l$ such that $k'_l >0$. Set $\mf{h}':=\mf{k}'-\mf{e}_l$, where
$\mf{e}_i:=(k_1,\ldots,k_n)$ and $k_j:=\delta^i_j$. Further, we set
$u_j:=(k'_j +1 -\delta_j^l)u_{\mf{k}'+\mf{e}_j-\mf{e}_l}$,
$u_0:=u_{\mf{h}'}$ and
$v'=\sum_{j=1}^nu_jg_j+u_0$. Comparing the coefficients of $g^{\mf{h}'}$ in
the equality
\begin{equation}\nonumber
0=\adR(E_i)v\simeq \sum\nolimits_{\mf{k}}u_{\mf{k}}\adR(E_i)g^{\mf{k}}+
\sum\nolimits_{\mf{h}}\adR(E_i)u_{\mf{h}}g^{\mf{h}}
\end{equation}
by using (\ref{eq-sim}) we obtain
\begin{equation}\nonumber
(k'_i+1-\delta^i_l)u_{\mf{k}'+\mf{e}_i-\mf{e}_l}E_i+
\adR(E_i)u_{\mf{h}'}=0,
\end{equation}
which means that $\adR(E_i)v'=0$. Since
$|{\mf{k}'+\mf{e}_j-\mf{e}_l}|=|v|$, we have
$u_j \in\mc{F}(\tUqg)$. Moreover, we also have
$\adR(E_i)u_j=0$ by construction and $\adR(E_i)g_j =\delta_i^j$.
Using these facts we get
$0=\adR(E_i)v' = u_iE_i + \adR(E_i)u_0$ and so
$\adR(E_i)u_0 \in \mc{F}(\tUqg)E_i \subseteq \mc{F}(\tUqg)\pchru$.
Recall that $\tUqg$ is the direct sum of $\adR$-invariant subspaces
$\mc{F}(\Uqg)\pchru f_\mu$. Therefore, it follows from
Corollary \ref{f-adREiinj} and the last relation that the nonvanishing
components of
$u_0 \in \tUqg$ are in $\mc{F}(\tUqg)\pchru$.
In particular, $u_0\tau(-2\lambda)\in\mc{F}(\tUqg)$ for
some $\lambda\in P_+(\mf{g})$. Hence there exists a natural number
$k$ such that $\adR(E_i^k)(u_0\tau(-2\lambda))=0$ and
$\adR(E_i^k)\tau(-2\lambda)=0$. Since $\adR(E_i)v'=0$, the latter implies
that $\adR(E_i^m)(v'\tau(-2\lambda))=v'\adR(E_i^m)\tau(-2\lambda)=0$ for all
$m\geq k$. Now we compute $\adR(E_i^m)(v'\tau(-2\lambda))$ by using the
equalities $\adR(E_i^m)(u_0\tau(-2\lambda))=0$ and
$\adR(E_i^k)\tau(-2\lambda)=0$. We then obtain
\begin{equation}\label{eq-szero}
\sum_{r=0}^{k-1}q_i^{-r(m-r)}\left[m\atop r\right]_{q_i}
u_i\prod_{j=1}^r(1-q_i^{-2j+2k})\tau(-2\lambda)E_i^r
\prod_{j=1}^{m-r-1}(1-q_i^{-2j})E_i^{m-r}=0.
\end{equation}
Next let us consider the equation
\begin{equation}\label{eq-sgen}
\sum_{r=s}^{k-1}q_i^{-r(m+s-r)}\left[m+s\atop r-s\right]_{q_i}
\prod_{j=1}^{m+s-r-1}(1-q_i^{-2j})\prod_{j=1}^r(1-q_i^{-2j+2k})\quad u_i=0
\end{equation}
for $s=0,1,\dots,k-1$. For $s=0$ equation (\ref{eq-sgen}) follows from
(\ref{eq-szero}),
because the algebra $\tUqg$ has no zero divisors. For general $s$,
(\ref{eq-sgen}) is easily
proved by induction. From (\ref{eq-sgen}) applied in the case $s=k-1$
we conclude that $u_i=0$,
because $q_i$ is not a root of unity.
Since $u_l:=k'_l u_{\mf{k}'}\neq 0$ by construction,
we arrived at a contradiction.
\end{bew}

\begin{satz}\label{s-FRGq}
If $q$ is transcendental, then we have
$\mc{F}(\RGqn)=\mc{F}(\Uqg)C$.
\end{satz}

\begin{bew}{}
Since $q$ is not a root of unity, the right adjoint action $\adR$ of the
algebra $\Uqg$ is a direct sum of simple submodules. Thus we have
$\mc{F}(\RGqn)\subseteq \mc{F}(\Uqg)C$ by Lemma \ref{l-RGq}.
In order to prove the
converse inclusion it suffices to check that $\adR(f_\mu)$ and $\adR(g_i)$
leave each space $cW(\lambda)$ for $c\in C$ and $\lambda\in P_+(\mf{g})$
invariant. Using the commutation rules (\ref{eq-kommutfg}) between
$f_\mu$, $g_i$ and the generators of $\Uqg$ this is easily done.
\end{bew}

\begin{lemma}\label{l-coideal}
For any finite dimensional $\adR$-invariant right coideal $\X$ of $\tUqg$
there exist elements $c\in C$ and $\lambda\in P_+(\mf{g})$ such that
$c\tau(-2\lambda) \in \X$.
\end{lemma}

\begin{bew}{}
Let $\{\mb{E}_{\mb{i}}\}$
and $\{\mb{F}_{\mb{j}}\}$
be vector space bases of $\Uqnp$
and $\Uqnm$, respectively, consisting of monomials in the generators
$E_1,\ldots,E_n$ and $F_1,\ldots,F_n$, respectively. Let $|{\mb{i}}|$
resp.\ $|{\mb{j}}|$ denote the $\mathbb{Z}$-degree of the monomial
$\mb{E}_{\mb{i}}$ resp.\ $\mb{F}_{\mb{j}}$. By Proposition \ref{s-hopfdual},
the set $\{\mb{E}_{\mb{i}}\mb{F}_{\mb{j}}f_\mu \}$
is a vector space basis of $\tUqg$. Fix a nonzero element $u$ of $\X$.
Then $u$ is a sum
$\sum_{{\mb{i}},{\mb{j}},\mu}\alpha_{{\mb{i}}{\mb{j}}\mu}
\mb{E}_{\mb{i}}\mb{F}_{\mb{j}}f_\mu$, where
$\alpha_{{\mb{i}}{\mb{j}}\mu}\in\mathbb{C}$. We choose indices
${\mb{i}}_0$ and
${\mb{j}}_0$ such that $\alpha_{{\mb{i}}_0
{\mb{j}}_0\mu_0}\not=0$ for some
$\mu_0\in(\mathbb{C}^\times)^n$ and that $\alpha_{{\mb{i}}_0{\mb{j}}\mu}=0$
if $|{\mb{i}}|>|{\mb{i}}_0|$
and $\alpha_{{\mb{i}}_0{\mb{j}}_0\mu}=0$ if $|{\mb{j}}|>|{\mb{j}}_0|$.
Since $\X$ is a right coideal,
$\Delta(u) \in \X \otimes \tUqg$. Hence the first tensor factor $\Delta(u)$
corresponding to the basis element
$\mb{E}_{{\mb{i}}_0}\mb{F}_{{\mb{j}}_0}f_{\mu_0}$ in the second tensor
factor belongs to $\X$. From the choice of the indices
${\mb{i}}_0,{\mb{j}}_0,\mu_0$ and the
formulas for the comultiplication of the generators $E_i,F_i,f_\mu$
it follows that this expression is of the form
$\alpha_{{\mb{i}}_0{\mb{j}}_0\mu_0}Kf_{\mu_0}$, where $K\in \chru$. Thus,
$f_\mu :=Kf_{\mu_0} \in \X \subset\mc{F}(\tUqg)$. By
Lemma \ref{l-adREik} applied with
$u=1$ and $\lambda =0$, there are $c\in C$ and $\gamma \in P(\mf{g})$ such that
$f_\mu = c\tau(-2\gamma)$.  Because $f_\mu =
c\tau(-2\gamma)\in\mc{F}(\tUqg)$, we have $\gamma \in P_+(\mf{g})$.
\end{bew}

The main result of the section is the following

\begin{thm}\label{t-adinvkoid}
Suppose that $q$ is transcendental. Then
any finite sum of the subspaces $cW(\lambda)$, where $c\in C$ and
$\lambda\in P_+(\mf{g})$, is a finite dimensional
$\adR$-invariant right coideal of $\RGqn$.
Conversely, each finite dimensional $\adR$-invariant right coideal $\X$ of
$\RGqn$ is a finite direct sum of subspaces
$cW(\lambda)$.
\end{thm}

\begin{bew}{}
{}In the proof of Proposition \ref{s-FRGq} we have already noted that
each space $cW(\lambda)$ is an $\adR$-invariant subspace of $\RGqn$. Since
\begin{displaymath}
\Delta(\adR(X)(c\tau(-2\lambda)))=
\adR(X_{(2)})(c\tau(-2\lambda))\otimes
S(X_{(1)})c\tau(-2\lambda)X_{(3)}
\end{displaymath}
for any $X\in\RGqn$, $cW(\lambda)$ is also a right coideal
of $\RGqn$.

Now let $\X$ be a finite dimensional $\adR$-invariant right coideal of
$\RGqn$. From Propositions \ref{s-locfinpart} and \ref{s-FRGq}
it follows that
there is a subset $\T$ of the product set
$C \times \pchru^{-1}$ such that
$\X\subset\bigoplus_{(c,\lambda)\in \T}\, cW({\lambda})\subset\tUqg$.
By Lemma \ref{l-coideal}, there are
elements $c_0 \in C$ and $\lambda_0 \in P_+(\mf{g})$ such that
$c_0\tau(-2\lambda_0) \in \X$. Consequently, the whole vector space
$c_0W(\lambda_0)$ is contained in $\X$. By the direct sum decomposition of
$\mc{F}(\RGqn)$ there is a finite dimensional $\adR$-invariant
subspace $\X'$ of $\bigoplus_{(c,\lambda)\in \T'}\, cW(\lambda)$
such that $\X=c_0W(\lambda_0)\oplus\X'$, where
$\T' := \T \setminus (c_0,\lambda_0)$.
As shown in the preceding paragraph each space $cW(\lambda)$ is a right
coideal. Hence $\bigoplus_{(c,\lambda)\in \T'}\, cW({\lambda})$
is a right coideal which has trivial intersection with $c_0W(\lambda_0)$.
This implies that $\X'$ is also a right coideal of $\RGqn$. Applying the
same reasoning with $\X$ replaced by $\X'$, the assertion follows by induction.
\end{bew}

An consequence of the preceding proof is the following

\begin{folg}\label{f-irredadinvkoid}
If $q$ is transcendental, then the irreducible $\adR$-invariant
right coideals of $\RGqn$
are precisely the subspaces $cW(\lambda)$, where
$\lambda\in P_+(\mf{g})$ and $c\in C$.
\end{folg}

Recall that by definition $\Uqg$ separates the elements of $\RGq$. Hence
$\RGqn$ separates the elements of $\RGq$. Thus, by
Proposition \ref{s-qliealg1}, there is a one-to-one correspondence between
finite dimensional $\adR$-invariant right coideals of
$\RGqn$ and finite dimensional
bicovariant differential calculi on $\RGq$. Therefore, the assertion of
Theorem \ref{t-adinvkoid}
gives a classification of all bicovariant differential calculi on the Hopf algebra
$\RGq$ by describing all possible quantum Lie algebras of such calculi. In
particular, by Corollary \ref{f-irredadinvkoid}, the subspaces
$cW(\lambda)$ are precisely the quantum Lie algebras of irreducible
bicovariant differential calculi on $\RGq$.

\section{Classification of bicovariant differential calculi on $\OGq$}
\label{sec-klassifikation}

Throughout this section $\OGq$ denotes the coordinate Hopf algebra
for one of the quantum groups $G_q=SL_q(n+1)$ or $G_q=Sp_q(2n)$ as defined
in \cite{a-FadResTak1} (see also \cite{b-KS}, Chap. 9)
and $\mf{g}$
is the corresponding Lie algebra $\liesl_{n+1}$ or $\liesp_{2n}$.
Further, we set $N=n+1$ if $G_q=SL_q(n+1)$ and
$N=2n$ if $G_q=Sp_q(2n)$.

\subsection{The main results}
\label{ss-mainresults}

Let us say that two pairs $(\zeta,v)$ and $(\tilde{\zeta},\tilde{v})$
of admissible $\zeta$ and $\tilde{\zeta}$ and of matrix
corepresentations $v$ and $\tilde{v}$ of $\OGq$ are
\textit{equivalent} if $\zeta=\tilde{\zeta}$ and $v$ and $\tilde{v}$
are equivalent. The main result of this paper is

\begin{thm}\label{t-main}
Let $\OGq$ be the coordinate Hopf algebra for one of the quantum groups
$SL_q(n+1)$ or $Sp_q(2n)$. Suppose that $q$ is a transcendental
complex number.
Then each finite dimensional bicovariant FODC $\Gamma$ over $\OGq$ is
isomorphic to a direct sum of bicovariant FODC $\Gamma_{\zeta_1}(v_1),
\ldots,\Gamma_{\zeta_k}(v_k)$ (see \ref{ss-OGqFODC}), where
$(v_i,\zeta_i)$ are mutually inequivalent pairs of admissible numbers
$\zeta_i$ for $G_q$ and finite dimensional irreducible matrix
corepresentations $v_i$ of $\OGq$ such that $(v_i,\zeta_i)\not=(1,1)$.
The number $k$ and the pairs
$(v_i,\zeta_i)$ are uniquely determined up to permutation of indices
and we have
\begin{displaymath}
\dim\Gamma=\sum\nolimits_{i=1}^k(\dim v_k)^2.
\end{displaymath}
The quantum Lie algebra $\X_\Gamma$ of $\Gamma$ is of the form
(\ref{eq-zentrgener}) with $c=c_{\zeta_1}(v_1)+\cdots+c_{\zeta_k}(v_k)$,
where $c_{\zeta_i}(v_i)$ is the central element of $\OGqn$ defined
by (\ref{eq-zentrinGamma}).
\end{thm}

Retaining the assumptions of Theorem \ref{t-main} we state three
corollaries.

\begin{folg} \label{f-FODCinner}
Any finite dimensional bicovariant FODC of $\OGq$ is inner.
\end{folg}

\begin{folg} \label{f-irredFODC}
For any irreducible finite dimensional bicovariant FODC $\Gamma$ there
exists an irreducible matrix corepresentation $v$ of $\OGq$ and an
admissible number $\zeta$ such that $\Gamma$ is isomorphic to
$\Gamma_\zeta(v)$. In particular, we then have $\dim\Gamma=(\dim v)^2$.
\end{folg}

\begin{folg} \label{f-FODCirred}
Let $v$ be a corepresentation of $\OGq$ and $\zeta$ an admissible
number. Suppose that $v$ is the direct sum
$\bigoplus_{i=0}^kn_iv_i$ of multiples of pairwise inequivalent
irreducible corepresentations $v_0,\ldots,v_k$, where $v_0$ denotes
the trivial corepresentation and $n_1\geq1,\ldots,n_k\geq1$.
If $\zeta\not=1$ and $n_0\geq1$, then the FODC $\Gamma_\zeta(v)$ is the
direct sum of the $k+1$ irreducible FODC
$\Gamma_\zeta(v_0),\Gamma_\zeta(v_1),\ldots,\Gamma_\zeta(v_k)$.
Otherwise, $\Gamma_\zeta(v)$ is the
direct sum of the $k$ irreducible FODC
$\Gamma_\zeta(v_1),\ldots,\Gamma_\zeta(v_k)$.
\end{folg}

\begin{bem}
The assertions of Proposition \ref{s-tensorFODC} and
Corollary \ref{f-FODCirred} are also valid for the quantum groups
$G_q=O_q(N)$.
\end{bem}

\subsection{Proof of the main results}
\label{ss-proofmain}
It is well known (\cite{b-KS}, Section 9.4) that there exists a unique
dual pairing $\langle \cdot ,\cdot \rangle$ of the Hopf algebras $\Uqg$
and $\OGq$
such that $\langle {l^+}^i_j,u^k_l\rangle =R^{ki}_{lj}$ and
$\langle {l^-}^i_j,u^k_l\rangle =R^{-1}{^{ik}_{jl}}$ for
$i,j,k,l=1,\ldots,N$, where ${l^\pm}^i_j$ are
the functionals (\ref{eq-generlfunc}) in case $u^i_j=v^i_j$ and $R$ is
the $R$-matrix for the vector representation of $\Uqg$.
Since $q$ is transcendental, the pairing
$\langle \cdot,\cdot\rangle $
is non-degenerate (\cite{b-KS}, Corollary 11.23). By the
definition of this pairing this implies that
the Hopf algebras $\OGq$ and $\RGq$ are isomorphic.

From the formulas for the comultiplication in $\Uqg$ and the
relations (\ref{eq-kommutfg}) it follows easily that the set
$C$ of central characters in $\RGqn \cong \OGqn$
consists of all functionals $f_\mu$, $\mu\in(\comp^\times)^n$, for which
${\tilde \mu_i}=1$
for all $i$. These are precisely the functionals
$\varepsilon_\zeta$ with admissible $\zeta$.

In the proof given below  we shall use some facts from the
corepresentation theory of the Hopf algebra
$\OGq$ (see, for instance, \cite{b-KS}).
One of the main results needed in the sequel is that there is a
one-to-one correspondence $v\to\varphi_v$ between irreducible
corepresentations $v$ of $\OGq$
and irreducible representations $\varphi_v$ of $\Uqg_{DJ}$
with dominant integral weights
$\lambda(v)\in P_+(\mf{g})$ (\cite{b-KS}, Theorem 11.22 and
Proposition 7.20). Using this correspondence the decomposition of
tensor products of irreducible corepresentations into direct sums of
irreducible components are given by the same formulas as in the classical
case. They are described by the graphical method in terms of the
Young frames (\cite{b-BarutRaczka}).

\begin{lemma}\label{l-coupling}
Let $v$ be an irreducible corepresentation of $\OGq$ with matrix
elements $v^i_j$, $i,j=1,\ldots,r$, with respect to some basis
of the underlying vector space. Suppose
that there exists an element $a(v)\in\mc{C}(v)$ such that
$l(a(v))=\tau(-2\lambda(v))$. Then for each admissible number $\zeta$
we have:\\
(i) $\X\cont(v)=W(\lambda(v))$ and the set
$\{l(v^i_j)\,|\,i,j=1,\ldots,r\}$ is a basis of the vector space
$\X\cont(v)$.\\
(ii) $\X_\zeta(v)+\mathbb{C}\varepsilon=
\varepsilon_\zeta W(\lambda(v))+\mathbb{C}\varepsilon$ and if
$(v,\zeta)\not=(1,1)$, then the set
$\{\varepsilon_\zeta l(v^i_j),\varepsilon\,|\,i,j=1,\ldots,r\}$
is a basis of the vector space $\X_\zeta(v)+\mathbb{C}\varepsilon$.
\end{lemma}

\begin{bew}{}
First we prove (i).
By Corollary \ref{f-irredadinvkoid}, $W(\lambda(v))$
is an irreducible $\adR$-invariant right coideal of $\RGqn \equiv\OGqn$.
Since
$l(a(v))=\tau(-2\lambda(v))\in \X\cont(v)\cap W(\lambda(v))$ by assumption
and $\X\cont(v)\cap W(\lambda(v))$ is an $\adR$-invariant right coideal
contained in $W(\lambda(v))$, it follows that $W(\lambda(v))\subset
\X\cont(v)$. Since
$\mathrm{dim}\,W(\lambda(v))=
r^2$ by
Proposition \ref{s-locfinpart} and $\mathrm{dim}\,\X\cont(v)\leq r^2$, the
preceding implies that $\X\cont(v)=W(\lambda(v))$ and that the
generating set $\{l(v^i_j)\}$ of $\X\cont(v)$ is a vector space basis of
$\X\cont(v)$.

Now we turn to (ii).
By (\ref{eq-ogqliealgfunc}), we have
$\X_\zeta(v)+\mathbb{C}\varepsilon=\varepsilon_\zeta\X\cont(v)+
\mathbb{C}\varepsilon$. Since $\X\cont(v)=W(\lambda(v))$
by (i), the first assertion follows.
If $(v,\zeta)\not=(1,1)$, then the sum
$\varepsilon_\zeta W(\lambda(v))+\mathbb{C}\varepsilon$
is direct by Theorem \ref{t-adinvkoid}.
Hence $\dim(\X_\zeta(v)+\mathbb{C}\varepsilon)=
\dim(\varepsilon_\zeta W(\lambda(v))+\mathbb{C}\varepsilon)=
r^2+1$. Therefore, the generating set $\{l(v^i_j),\varepsilon\}$ of
$\X_\zeta(v)+\mathbb{C}\varepsilon$ consisting of $r^2+1$
functionals is a vector space basis of
$\X_\zeta(v)+\mathbb{C}\varepsilon$.
\end{bew}

In the following example we shall show that for the irreducible
corepresentation $v_k$
with Young frame $[1^k,0^{n-k}]$,
$k=1,\ldots,n,$ there is an element $a(v_k)\in\mc{C}(v_k)$ satisfying
$l(a(v_k))=\tau(-2\lambda(v_k))$.

\begin{bsp}\label{b-antisymm}
For $G_q=SL_q(N)$ and $G_q=Sp_q(N)$, respectively, let $V_q$ denote the
quantum vector space $\mathbb{C}_q^N$ and the quantum symplectic
space $Sp_q^N$, respectively, and $\Lambda(V_q)$ their exterior
algebras (see \cite{a-FadResTak1} or \cite{b-KS}, Definitions 9.4 and
9.12). Let $y_1,\ldots,y_N$ be the generators of $\Lambda(V_q)$ and
$\Lambda(V_q)_k$ the subspace of $\Lambda(V_q)$ of
homogeneous elements $y\in\Lambda(V_q)$ of degree $k$.
Recall that the right coaction $\varphi_R$ of $\OGq$ on
$\Lambda(V_q)$ is the algebra homomorphism
$\varphi_R:\Lambda(V_q)\to\Lambda(V_q)\otimes\OGq$
determined by $\varphi_R(y_i)=y_j\otimes u^j_i$.
Let $\varphi_{R,k}$ denote the restriction of $\varphi_R$ to
the invariant subspace $\Lambda(V_q)_k$.
Since the set $\{y_I\,|\,I=(i_1,\ldots,i_k),i_1<i_2<\cdots<i_k\}$
forms a basis in $\Lambda(V_q)_k$, there exist elements
$\mc{D}^I_J\in\OGq$, such that
$\varphi_{R,k}(y_I)=y_J\otimes\mc{D}^J_I$.
Then we have
\begin{equation}\label{eq-deltaDIJ}
\Delta(\mc{D}^I_J)=\mc{D}^I_M\otimes\mc{D}^M_J,
\end{equation}
where the summation is over all sets
$M=(m_1,\ldots,m_k)$ of integers such that
$1 \leq m_1<m_2<\cdots<m_k \leq N$. In the case $G_q=SL_q(N)$, this
is proved in \cite{b-KS}, 9.2.2, but all
considerations remain valid for $G_q=Sp_q(N)$ as well.

Now we suppose that $k\leq n$. Let us take a closer look at the
elements $\mc{D}^I_J$.
If $i_k\leq n$, then
\begin{equation}
\mc{D}^I_J=\sum_{\sigma\in\mc{P}_k}
(-q)^{\ell(\sigma)}u^{\sigma(i_1)}_{j_1}\cdots u^{\sigma(i_k)}_{j_k},
\end{equation}
where $\mc{P}_k$ is the group of permutations of the set
$\{i_1,\ldots,i_k\}$ and $\ell(\sigma)$ is the length of the permutation
$\sigma$. For $G_q=SL_q(n+1)$ the latter is proved in \cite{b-KS},
Proposition 9.7. For $G_q=Sp_q(2n)$ the proof is similar (it suffices
to use the first two equations in \cite{b-KS}, Proposition 9.17(i)).
Let $\mc{O}_>(G_q)$ denote the non-unital subalgebra of $\OGq$
generated by the matrix elements $u^i_j$, $i>n$, $j\leq n$.
We easily see that
if $i_k>n$ and $j_k\leq n$, then $\mc{D}^I_J\in\mc{O}_>(G_q)$.

Now let $v_k$ be the irreducible corepresentation
with Young frame $[1^k,0^{n-k}]$, $1\leq k\leq n$. Set
\begin{equation}
\mc{D}_{q,k}:=\mc{D}^I_J=\sum_{\sigma\in\mc{P}_k}(-q)^{\ell(\sigma)}
u^{\sigma(1)}_1\cdots u^{\sigma(k)}_k,
\end{equation}
where $I=J=(1,\ldots,k)$.
Obviously, $\mc{D}_{q,k} \in \mc{C}(u^{\otimes k})$. The Drinfeld-Jimbo algebra $\Uqg_{DJ}$
acts on $\OGq$ by $\varphi(f)a=a_{(1)}\langle f,a_{(2)}\rangle $, $a\in\OGq$,
$f\in\Uqg_{DJ}$.
{}From the formula
$\Delta(F_i)=F_i\otimes 1+K_{\alpha_i}^{-1}\otimes F_i$
for the comultiplication of $\Uqg_{DJ}$
we obtain
\begin{equation}
\Delta^{(k-1)}(F_i)=\sum_{j=1}^k{K_{\alpha_i}^{-1}}^{\otimes j-1}
\otimes F_i\otimes 1^{\otimes k-j},\quad k\geq2.
\end{equation}
Recall that $\Delta(\mc{D}_{q,k})$ is given by
(\ref{eq-deltaDIJ}).
Since $\langle F_i,u^r_s\rangle =\delta^r_i\delta^s_{i+1}$ for $s\leq n$ and
$\langle K_{\alpha_i},u^r_s\rangle =0$ for $r\not=s$, this implies
that $\varphi(F_i)\mc{D}^M_J=0$ for all $M=(m_1,\ldots,m_k)$ and
$i=1,\ldots,n$. Hence we have $\varphi(F_i)\mc{D}_{q,k}=0$ for
$i=1,\ldots,n$.
Using the explicit formula
$\langle K_{\alpha_i},u^r_r\rangle =q_i^{-\delta^i_r+\delta^{i+1}_r}$ for
$r\leq n$ (see \cite{b-KS}, Section 9.4)
one easily verifies that
$\varphi(K_{\alpha_i}^{-1})\mc{D}_{q,k}=q_i^{\delta_{ik}}\mc{D}_{q,k}=
q^{(\alpha_i,\fuw_k)}\mc{D}_{q,k}$.
That is, $\mc{D}_{q,k}$ is a highest weight vector with highest
weight $\fuw_k$ with respect to the ordered sequence of simple
roots $-\alpha_1,\ldots,-\alpha_n$ for the representation $\varphi$
of $\Uqg_{DJ}$. Therefore,
by the Peter-Weyl decomposition of the coordinate Hopf algebra $\OGq$
(\cite{b-KS}, Theorem 11.22), $\mc{D}_{q,k}$ belongs to the
coalgebra $\mc{C}(v_k)$.

Suppose that $I=J=\{1,\ldots,k\}$. {}From (\ref{eq-deltaDIJ})
we get $l(\mc{D}^I_J)=\sum\nolimits_M S(l^-(\mc{D}^I_M))l^+(\mc{D}^M_J)$.
If $M$ is a multi-index such that $m_k>n$, then we have
$\mc{D}^M_J\in\mc{O}_>(G_q)$ and hence $l^+(D^M_J)=0$.
Therefore, it suffices to sum over
multi-indices $M=(m_1,\ldots,m_k)$ with $m_k\leq n$. Further, from the
explicit formulas (see for instance \cite{b-KS}, 8.5.4) we know that
$K_i={l^-}^1_1\cdots{l^-}^i_i, i=1,\ldots,n$.
Now we compute
\begin{align*}
l(\mc{D}^I_J)&=\sum_{m_1<m_2<\ldots<m_k\leq n}
\sum_{\sigma,\sigma'}(-q)^{\ell(\sigma)+\ell(\sigma')}
S(l^-(u^{\sigma(1)}_{m_1}\cdots u^{\sigma(k)}_{m_k}))
l^+(u^{\sigma'(m_1)}_1\cdots u^{\sigma'(m_k)}_k)\\
&=\sum_{m_1<m_2<\ldots<m_k\leq n}
\sum_{\sigma,\sigma'}(-q)^{\ell(\sigma)+\ell(\sigma')}
S({l^-}^{\sigma(k)}_{m_k}\cdots {l^-}^{\sigma(1)}_{m_1})
{l^+}^{\sigma'(m_k)}_k\cdots {l^+}^{\sigma'(m_1)}_1\\
&=S({l^-}^k_k\cdots{l^-}^1_1){l^+}^k_k\cdots{l^+}^1_1=
({l^+}^1_1\cdots{l^+}^k_k)^2=\tau(-2\fuw_k)=\tau(-2\lambda(v_k)).
\end{align*}

Moreover, the last reasoning shows also that
\begin{equation}\label{eq-Slmtensorlpvonav}
S(l^-(a(v)_{(1)}))\otimes l^+(a(v)_{(2)})={l^+}^1_1\cdots{l^+}^k_k\otimes
{l^+}^1_1\cdots{l^+}^k_k.
\end{equation}
\end{bsp}

The crucial step in the proof of Theorem \ref{t-main} is the following
lemma.

\begin{lemma}\label{l-findeav}
For each irreducible corepresentation $v$ of $\OGq$ there exists an
element $a(v)\in\mc{C}(v)$ such that $l(a(v))=\tau(-2\lambda(v))$.
\end{lemma}

\begin{bew}{}
For $v=1$ we obtain $l(1)=\tau(0)=\varepsilon$. If $v=v_k$ is the
irreducible corepresentation corresponding to the Young frame
$[1^k,0^{n-k}]$, $k=1,\ldots,n$, then we have shown in
Example \ref{b-antisymm}
that $a(v_k)=\mc{D}_{q,k}$ is an element of $\mc{C}(v)$ such
that $l(a(v_k))=\tau(-2\lambda(v_k))=\tau(-2\fuw_k)$.

The general case will be treated by an induction procedure with respect
to an ordering of the Young frames. In order to do so we first introduce
some notation.

Let us consider the Young frame of an irreducible corepresentation $v$ of
$\OGq$. Let $m_j(v)$ denote the number of columns of length $j$.
Then $m(v):=\sum_jm_j(v)$ is the number of columns of the Young frame
of $v$ and $\lambda(v):=\sum_jm_j(v)\fuw_j$ is the highest weight
corresponding to $v$ with respect to the simple roots
$-\alpha_1,\ldots,-\alpha_n$.
Let $k(v)$ be the number of squares in the last (i.\,e.\ $m(v)$-th)
column. Let $\succ$ denote the lexicographic ordering of pairs $(m,k)$,
where $m\in\mathbb{N}$, $k\in\{1,\ldots,n\}$. That is,
$(m,k)\succ(m',k')$ if either $m>m'$ or $m=m'$, $k>k'$.
As shown in the first part of this proof, the assertion is true for all
irreducible corepresentations $v$ such that $m(v)\leq 1$.
Now let $v$ be an arbitrary irreducible corepresentation of $\OGq$
such that $(m(v),k(v))\succ(1,n)$. Suppose that the assertion holds for all
irreducible corepresentations $w$ such that
$(m(v),k(v))\succ(m(w),k(w))$. We shall prove the assertion for $v$, that
is, we have to show that $\tau(-2\lambda(v))\in\X\cont(v)$.

Let $v_1$ be the irreducible corepresentation with Young frame consisting
of the last column of $v$ (i.\,e.\ with Young frame $[1^{k(v)},0^{n-k(v)}]$)
and let $w$ be the irreducible corepresentation which is obtained if the
last column in the Young frame of $v$ is cancelled.
Then we have $\lambda(v)=\lambda(w)+\lambda(v_1)$ and
$\lambda(v_1)=\fuw_{k(v)}$. Since $(m(v),k(v))\succ(m(w),k(w))$,
by assumption there is an element $a(w)\in\mc{C}(w)$ such that
$l(a(w))=\tau(-2\lambda(w))$. For $a,b\in\OGq$, we easily compute that
$l(ab)=S(l^-(a_{(1)}))l(b)l^+(a_{(2)})$.
Setting $a=a(v_1)=\mc{D}_{q,k(v)}$ and $b=a(w)$ and using
(\ref{eq-Slmtensorlpvonav}) we get
\begin{equation}\label{eq-lvonavaw}
\begin{split}
l(a(v_1)a(w))&={l^+}^1_1\cdots{l^+}^{k(v)}_{k(v)}\tau(-2\lambda(w))
{l^+}^1_1\cdots{l^+}^{k(v)}_{k(v)}\\
&=\tau(-2\lambda(w))({l^+}^1_1\cdots{l^+}^{k(v)}_{k(v)})^2=
\tau(-2\lambda(w)-2\fuw_{k(v)})=\tau(-2\lambda(v)).
\end{split}
\end{equation}
Since the Hopf algebra $\OGq$ is cosemisimple (\cite{a-Hayashi}, or
\cite{b-KS}, Theorem 11.22), the corepresentation
$v_1\otimes w$ of $\OGq$ decomposes into a direct sum $\bigoplus_jw_j$
of irreducible corepresentations $w_j$. {}From the decomposition rules
of the tensor product in terms of Young frames (see, for instance,
\cite{b-BarutRaczka}, \S8,c) we see that precisely one summand,
say $w_0$, has the same Young
frame as $v$ and that $(m(v),k(v))\succ(m(w_j),k(w_j))$ for all other
summands. Thus, $w_0$ is equivalent to $v$, so that
$\mc{C}(v)=\mc{C}(w_0)$. By the induction hypothesis combined with
Lemma \ref{l-coupling}(i) we conclude that
$\X\cont(w_j)=W(\lambda(w_j)),j\not=0$.
The element $a(v_1)a(w)$ is in $\mc{C}(v_1\otimes w)$ by
construction and hence
$\tau(-2\lambda(v))=l(a(v_1)a(w))\in\X\cont(v_1\otimes w)$
by (\ref{eq-lvonavaw}). Since $W(\lambda(v))$ is the irreducible
$\adR$-invariant right coideal generated by $\tau(-2\lambda(v))$ and
$\X\cont(v_1\otimes w)$ is also an $\adR$-invariant
right coideal, the latter implies that
$W(\lambda(v))\subset\X\cont(v_1\otimes w)$.
{}From the preceding
and the facts that $v_1\otimes w=\bigoplus_jw_j$ and
$\mc{C}(v)=\mc{C}(w_0)$ we obtain
\begin{equation} \label{eq-tensorsplitliealg}
W(\lambda(v))\subset\X\cont(v)+\sum_{j\not=0}W(\lambda(w_j)).
\end{equation}
On the other hand, $\X\cont(v)$ is an $\adR$-invariant
right coideal of $\Uqg$
and hence a sum $\sum_iW(\lambda_i)$ of irreducible $\adR$-invariant
right coideals $W(\lambda_i)$, $\lambda_i\in P_+(\mf{g})$.
Therefore, since the sum
$\sum_{\lambda\in P_+(\mf{g})}W(\lambda)$ is a \textit{direct sum}
by Proposition \ref{s-locfinpart}, equation (\ref{eq-tensorsplitliealg})
implies that $W(\lambda(v))\subset\X\cont(v)$.
Hence we have $\tau(-2\lambda(v))\in W(\lambda(v))\subset\X\cont(v)$,
that is, $\tau(-2\lambda(v))=l(a(v))$ for some element $a(v)\in\mc{C}(v)$.
\end{bew}

\begin{folg}\label{f-unabhfunc}
Let $v=(v^i_j)_{i,j=1,\ldots,m}$ be an irreducible corepresentation
of $\OGq$ and let $\zeta\in\mathbb{C}$ be admissible for $G_q$
such that $(v,\zeta)\not=(1,1)$. Then the linear functionals
$X_{ij}=\varepsilon_\zeta l(v^i_j)-\delta_{ij}\varepsilon$,
$i,j=1,\ldots,m$, are linearly independent. Moreover, the sets
of linear functionals
$\{\varepsilon_\zeta l(v^i_j),\varepsilon\,|\,i,j=1,\ldots,m\}$ and
$\{l(v^i_j)\,|\,i,j=1,\ldots,m\}$ are also linearly independent.
\end{folg}

\begin{bew}{}
The second assertion follows from Lemma \ref{l-findeav} and
Lemma \ref{l-coupling}(ii) and (i).
The first one is a consequence of the linear independence of
functionals $\varepsilon_\zeta l(v^i_j),\varepsilon$.
\end{bew}

\begin{lemma}\label{l-zentrnichtnull}
Let $v=(v^i_j)_{i,j=1,\ldots,m}$ be an irreducible corepresentation
of $\OGq$ and let $\zeta\in\mathbb{C}$ be admissible for $G_q$
such that $(v,\zeta)\not=(1,1)$. Then the central functionals
$c_\zeta(v)-c_\zeta(v)(1)\varepsilon\in\X_\zeta(v)$ defined by
(\ref{eq-zentrinGamma}) are nonzero.
\end{lemma}

\begin{bew}{}
By (\ref{eq-zentrinGamma}) the element
$c_\zeta(v)-c_\zeta(v)(1)\varepsilon$ is a linear combination of
functionals $\varepsilon_\zeta l(v^i_j)$, $i,j=1,\ldots,m$, and
$\varepsilon$. By Corollary \ref{f-unabhfunc}, these functionals are
linearly independent. Therefore, since the matrix
$D^{-1}$ is invertible,
$c_\zeta(v)-c_\zeta(v)(1)\varepsilon\not=0$.
\end{bew}

Now we are able to give the proofs of our main results.

\begin{bew}{ of Theorem \ref{t-main}}
Let $\Gamma$ be a finite dimensional bicovariant FODC over $\OGq$
and $\X_\Gamma$ be its quantum Lie algebra. By Proposition
\ref{s-qliealg1}(i), $\X_\Gamma+\mathbb{C}\varepsilon$ is an
$\adR$-invariant right coideal of the Hopf dual $\OGqn\equiv \RGqn$.
Therefore, by Theorem \ref{t-adinvkoid},  $\X_\Gamma+
\mathbb{C}\varepsilon$ is a direct sum $\bigoplus_{i=0}^k
\varepsilon_{\zeta_i}W(\lambda_i)$, where $\lambda_i\in P_+(\mf{g})$
and $\zeta_i$ is an admissible number for $G_q$. Without loss of
generality we can assume that $\mathbb{C}\varepsilon=
\varepsilon_{\zeta_0}W(\lambda_0)$, that is, $\zeta_0=1$ and
$\lambda_0=0$. Let $v_i$ denote the irreducible corepresentation
of $\OGq$ which corresponds to the irreducible representation of
$\Uqg$ with highest weight $\lambda_i$ (with respect to the
ordered sequence $\{-\alpha_1,\ldots,-\alpha_n\}$ of simple roots).
Since $\bigoplus_{i=0}^k\varepsilon_{\zeta_i}W(\lambda_i)$ is a direct
sum, it follows that the pairs $(\zeta_i,v_i)$, $i=0,\ldots,k$,
are mutually inequivalent. {}From Lemma \ref{l-coupling}(ii) and
Lemma \ref{l-findeav} we conclude that $\X_{\zeta_i}(v_i)+
\mathbb{C}\varepsilon=\varepsilon_{\zeta_i}W(\lambda_i)+
\mathbb{C}\varepsilon$ for $i=1,\ldots,k$. Therefore, since
$\mathbb{C}\varepsilon\oplus\sum_{i=1}^k\X_{\zeta_i}(v_i)=
\bigoplus_{i=0}^k\varepsilon_{\zeta_i}W(\lambda_i)=
\mathbb{C}\varepsilon\oplus\X_\Gamma$, it follows that
$\sum_{i=1}^k\X_{\zeta_i}(v_i)$ is a direct sum and equal to
$\X_\Gamma$. Since $\X_{\zeta_i}(v_i)$ is the quantum Lie algebra
of the FODC $\Gamma_{\zeta_i}(v_i)$, this implies that the sum
of bicovariant FODC $\Gamma_{\zeta_i}(v_i)$, $i=1,\ldots,k$, is a
direct sum and that the quantum Lie algebras of the FODC
$\bigoplus_{i=1}^k\Gamma_{\zeta_i}(v_i)$ and $\Gamma$ coincide.
Hence $\Gamma$ is isomorphic to
$\bigoplus_{i=1}^k\Gamma_{\zeta_i}(v_i)$.

Next we prove that $\X_\Gamma=\X[c]$ with $c=c_{\zeta_1}(v_1)+
\cdots+c_{\zeta_k}(v_k)$. Since $c_i:=c_{\zeta_i}(v_i)-
c_{\zeta_i}(v_i)(1)\varepsilon\in\X_{\zeta_i}(v_i)$ as noted in
\ref{ss-OGqFODC}, we have $c-c(1)\varepsilon\in\X_\Gamma$
and so $c\in\X_\Gamma\oplus\mathbb{C}\varepsilon$. Let us write
$\Delta(c)=\sum_jx_j\otimes y_j$ with $\{y_j\}$ linearly independent.
Since $\X_\Gamma\oplus\mathbb{C}\varepsilon$ is a right coideal by
Proposition \ref{s-qliealg1}(i), $x_j\in\X_\Gamma\oplus
\mathbb{C}\varepsilon$ and hence $x_j-x_j(1)\varepsilon\in\X_\Gamma$.
As noted in \ref{ss-centrgenFODC}, $\X[c]$ is the linear span of
functionals $x_j-x_j(1)\varepsilon$, so we conclude that
$\X[c]\subset\X_\Gamma$ and hence
\begin{equation}\label{eq-centrgenercoid}
\X[c]\oplus\mathbb{C}\varepsilon\subset\X_\Gamma\oplus
\mathbb{C}\varepsilon=\bigoplus_{i=1}^k\X_{\zeta_i}(v_i)\oplus
\mathbb{C}\varepsilon=
\bigoplus_{i=0}^k\varepsilon_{\zeta_i}W(\lambda_i).
\end{equation}
On the other hand, from Proposition \ref{s-centrgenFODC} and
Proposition \ref{s-qliealg1}(i) we obtain that
$\X[c]\oplus\mathbb{C}\varepsilon$ is an $\adR$-invariant right
coideal of $\OGqn$ and so a direct sum of irreducible right coideals
$\varepsilon_\zeta W(\lambda)$. Equation (\ref{eq-centrgenercoid})
implies that $\X[c]\oplus\mathbb{C}\varepsilon=
\bigoplus_{i\in I}\varepsilon_{\zeta_i}W(\lambda_i)\oplus
\mathbb{C}\varepsilon$ for some subset $I$ of $\{1,\ldots,k\}$.
Setting $a=1$ in (\ref{eq-zentrgener}) we get $c-c(1)\varepsilon
\in\X[c]$ and hence $c=c_1+\cdots+c_k+c(1)\varepsilon\in\X(c)
\oplus\mathbb{C}\varepsilon=\bigoplus_{i\in I}\X_{\zeta_i}(v_i)
\oplus\mathbb{C}\varepsilon$. Since $c_i\not=0$ by Lemma
\ref{l-zentrnichtnull} and $c_i\in\X_{\zeta_i}(v_i)$ for
$i=1,\ldots,k$, the latter is only possible if $I=\{1,\ldots,k\}$.
This implies at once that $\X[c]=\X_\Gamma$ and completes the proof
of Theorem \ref{t-main}.
\end{bew}


\begin{bew}{ of Corollary \ref{f-FODCinner}}
Since the linear functionals (\ref{eq-ogqliealgfunc}) for any FODC
$\Gamma_{\zeta_i}(v_i)$, $i=1,\ldots,k$, are linearly independent
by Corollary \ref{f-unabhfunc}, we have $\Gamma_{\zeta_i}(v_i)=
\tilde{\Gamma}$ (in the notation of \ref{ss-OGqFODC}, see
(\ref{eq-tildeGamma})). Hence the corresponding biinvariant element
$\theta_i$ belongs to $\Gamma_{\zeta_i}(v_i)$. By the definition
of the direct sum FODC, $\theta:=\sum_{i=1}^k\theta_i$ is a
biinvariant element of $\bigoplus_{i=1}^k\Gamma_{\zeta_i}(v_i)$ which
defines the differentiation of the FODC
$\bigoplus_{i=1}^k\Gamma_{\zeta_i}(v_i)\simeq\Gamma$ by
(\ref{eq-innerdif}). Thus $\Gamma$ is inner.
\end{bew}

\begin{bew}{ of Corollary \ref{f-irredFODC}}
Since the bicovariant FODC $\Gamma$ is irreducible,
by Lemma \ref{l-equivirredFODC} there is no $\adR$-invariant
right coideal $\mc{Y}$ of $\OGqn$ such that $\mathbb{C}\varepsilon
\subset\mc{Y}\subset\X_\Gamma\oplus\mathbb{C}\varepsilon$ and
$\mc{Y}\not=\mathbb{C}\varepsilon,
\X_\Gamma\oplus\mathbb{C}\varepsilon$. Therefore we have $k=1$
in the proof of Theorem \ref{t-main} and so $\Gamma$ is
isomorphic to $\Gamma_{\zeta_1}(v_1)$.
\end{bew}

\begin{bew}{ of Corollary \ref{f-FODCirred}}
The assumption $v=\bigoplus_{i=0}^kn_iv_i$ yields that
$\X_\zeta(v)=\sum_{i=0}^kX_\zeta(v_i)$ for $\zeta\not=1$ and $n_0>0$
and $\X_\zeta(v)=\sum_{i=1}^kX_\zeta(v_i)$ otherwise.
Since the corepresentation
$v_0,\ldots,v_k$ are pairwise inequivalent, $\lambda(v_i)\not=
\lambda(v_j)$ for $i\not=j$. Hence $\sum_{i=0}^kW(\lambda_i)$
($\sum_{i=0}^kW(\lambda_i)$, resp.) is a direct sum. This implies that
$\X_\zeta(v)=\bigoplus_{i=0}^k\X_\zeta(v_i)$
($\X_\zeta(v)=\bigoplus_{i=1}^k\X_\zeta(v_i)$, resp.).
By Lemma \ref{l-direktesumme}, $\Gamma_\zeta(v)$ is a direct sum of
$\Gamma_\zeta(v_0),\ldots,\Gamma_\zeta(v_k)$
($\Gamma_\zeta(v_1),\ldots,\Gamma_\zeta(v_k)$, resp.).
\end{bew}

\section{Factorizability of the Hopf algebras $\mc{O}(G_q)$}

Following \cite{b-Majid1}, we say that a pair $(\A,\br)$ of a
coquasitriangular Hopf algebra $\A$ and a universal $r$-form $\br$
is \textit{factorizable} if the bilinear form $\mb{q}$ on
$\A\otimes\A$ defined by
\begin{displaymath}
\mb{q}(a\otimes b)=\br(b_{(1)}\otimes a_{(1)})\br(a_{(2)}\otimes b_{(2)}),\quad
a,b\in\A,
\end{displaymath}
is non-degenerate. The corresponding dual notion for quasitriangular
Hopf algebras was introduced in \cite{a-ResSem}.

The functionals $l^i_j$ (see (\ref{eq-generlfunc}))
are of the form $l^i_j(\cdot)=\mb{q}(\cdot\otimes v^i_j)$, see
\cite{b-KS}. For arbitrary $a\in\A$, we define a
linear functional $l(a)$ on $\A$ by
$l(a)(\cdot):=\mb{q}(\cdot\otimes a)$.

\begin{thm}\label{s-factoris}
Suppose that $q$ is transcendental. Let $\br$ be the
canonical universal $r$-form on $\OGq$ determined by (\ref{eq-rform}).
Then the pair $(\OGq,\br)$, $G_q=SL_q(n+1),Sp_q(2n)$, is factorizable.
\end{thm}

\begin{bew}{}
Suppose that $l(a)=0$ for some $a\in\OGq$. Let $v(\lambda)$
denote the corespresentation of $\OGq$ corresponding to $\lambda\in
P_+(\mf{g})$. Because $\OGq$ is cosemisimple (see \cite{a-Hayashi} or
\cite{b-KS}), $a$ is a finite sum of elements
$a_\lambda\in\mc{C}(v(\lambda))$. Recall that
$l(a_\lambda)\in\X\cont(v(\lambda))$ by definition. Since
$\X\cont(v(\lambda))=W(\lambda(v))$ by Lemmas \ref{l-coupling} and
\ref{l-findeav} and the sum of spaces $W(\lambda)$
is direct, the assumption
$l(a)=\sum_\lambda l(a_\lambda)=0$ implies that $l(a_\lambda)=0$ for
all $\lambda$. Since the functionals $l(v(\lambda)^i_j)$ form a basis of
$\X\cont(v(\lambda))$ by Lemmas \ref{l-coupling} and
\ref{l-findeav}, $a_\lambda=0$ for all $\lambda$. Thus we have shown
that $l(a)(b)=\mb{q}(a\otimes b)=0$ for all
$b\in\OGq$ implies that $a=0$. Since $\mb{q}(S(b)\otimes
S(a))=\mb{q}(a\otimes b)$
and $S$ is invertible, $\mb{q}$ is non-degenerate.
\end{bew}

\nocite{a-Jurco1}\nocite{a-CSchWW1}


{\small
Address: Universit{\"a}t Leipzig, Fakult{\"a}t f{\"u}r Mathematik und
Informatik und NTZ, Augustusplatz 10/11, D--04109 Leipzig, Germany}

\end{document}